\RequirePackage{iftex}
\ifPDFTeX
  \documentclass[a4paper]{amsart}
\else
  \documentclass[dvipdfmx, a4paper]{amsart}
\fi

\usepackage{amsmath, amssymb}
\usepackage{mathtools}
\usepackage{tikz-cd}

\usepackage{comment}
\usepackage{xcolor}
\usepackage{enumitem}
\setlist[enumerate,1]{
    label=(\roman*),
    ref=(\roman*),
    font=\normalfont,
    topsep=3pt,
    itemsep=2pt,
    parsep=0pt,
    partopsep=0pt,
    leftmargin=*,
    widest=iii
}
\setlist[enumerate,2]{
    label=(\alph*),
    ref=\theenumi(\alph*),
    font=\normalfont,
    topsep=2pt,
    itemsep=1pt,
    parsep=0pt,
    leftmargin=*,
    widest=b
}
\usepackage{graphicx}
\usepackage{booktabs}

\usepackage[colorlinks=false, linkbordercolor=green, citebordercolor=red, hypertexnames=false]{hyperref}
\usepackage{bookmark}
\usepackage{amsthm}

\usepackage{orcidlink}

\usepackage[noabbrev, capitalize]{cleveref}

\makeatletter
\AtBeginEnvironment{fact}{\crefalias{theorem}{fact}}
\AtBeginEnvironment{lemma}{\crefalias{theorem}{lemma}}
\AtBeginEnvironment{proposition}{\crefalias{theorem}{proposition}}
\AtBeginEnvironment{corollary}{\crefalias{theorem}{corollary}}
\AtBeginEnvironment{definition}{\crefalias{theorem}{definition}}
\AtBeginEnvironment{example}{\crefalias{theorem}{example}}
\AtBeginEnvironment{remark}{\crefalias{theorem}{remark}}
\AtBeginEnvironment{question}{\crefalias{theorem}{question}}
\AtBeginEnvironment{conjecture}{\crefalias{theorem}{conjecture}}
\makeatother

\crefformat{equation}{(#2#1#3)}
\Crefformat{equation}{Equation (#2#1#3)}

\crefname{theorem}{Theorem}{Theorems}
\Crefname{theorem}{Theorem}{Theorems}
\crefname{lemma}{Lemma}{Lemmas}
\Crefname{lemma}{Lemma}{Lemmas}
\crefname{proposition}{Proposition}{Propositions}
\Crefname{proposition}{Proposition}{Propositions}
\crefname{corollary}{Corollary}{Corollaries}
\Crefname{corollary}{Corollary}{Corollaries}
\crefname{fact}{Fact}{Facts}
\Crefname{fact}{Fact}{Facts}
\crefname{definition}{Definition}{Definitions}
\Crefname{definition}{Definition}{Definitions}
\crefname{example}{Example}{Examples}
\Crefname{example}{Example}{Examples}
\crefname{remark}{Remark}{Remarks}
\Crefname{remark}{Remark}{Remarks}
\crefname{question}{Question}{Questions}
\Crefname{question}{Question}{Questions}
\crefname{conjecture}{Conjecture}{Conjectures}
\Crefname{conjecture}{Conjecture}{Conjectures}

\numberwithin{equation}{section}

\theoremstyle{definition}
\newtheorem{theorem}{Theorem}[section]
\newtheorem{lemma}[theorem]{Lemma}
\newtheorem{corollary}[theorem]{Corollary}
\newtheorem{proposition}[theorem]{Proposition}
\newtheorem{fact}[theorem]{Fact}
\newtheorem{definition}[theorem]{Definition}
\newtheorem{example}[theorem]{Example}
\newtheorem{remark}[theorem]{Remark}

\newtheorem{conjecture}[theorem]{Conjecture}
\newtheorem*{acknowledgements}{Acknowledgements}
\newtheorem*{aistatement}{Use of AI statement}

\makeatletter
\newtheorem*{rep@theorem}{\pointerep@title}
\newenvironment{reptheorem}[3][]{%
  \def\pointerep@title{#2 \ref{#3}}%
  \if\relax\detokenize{#1}\relax\begin{rep@theorem}\else\begin{rep@theorem}[#1]\fi}
  {\end{rep@theorem}}
\makeatother

\DeclareEmphSequence{\itshape}

\usepackage{marginnote}

\newcommand{\checkmarker}{\marginnote{\textcolor{magenta}{\scalebox{0.8}{\bfseries [Check!]}}}}

\usepackage{autonum}

\newcommand{\Z}{\mathbf{Z}}
\newcommand{\R}{\mathbf{R}}
\newcommand{\N}{\mathbf{N}}

\DeclareMathOperator{\diam}{diam}
\DeclareMathOperator{\pr}{pr}

\DeclarePairedDelimiter{\abs}{\lvert}{\rvert}
\DeclarePairedDelimiter{\norm}{\lVert}{\rVert}

\newcommand{\CM}{\mathcal{C}_M}
\newcommand{\dJ}{d_J}
\DeclareMathOperator{\dom}{dom}

\DeclareMathOperator{\Ric}{Ric}
\newcommand{\cdf}{\operatorname{cd}}

\csname endofdump\endcsname

\begin{document}

\title[Bartnik's conjecture and the space of Cauchy hypersurfaces]
{Bartnik's splitting conjecture and the coarse geometry of the space of Cauchy hypersurfaces}
\author[Keita Takahashi]{Keita Takahashi}
\address{RIKEN Center for Interdisciplinary Theoretical and Mathematical Sciences (iTHEMS), RIKEN, Wako 351-0198, Japan}
\address{Department of Mathematics, Institute of Science Tokyo,
  O-okayama, Meguro, Tokyo 152-8551, Japan}
\email{keita.takahashi@riken.jp}
\date{\today}

\begin{abstract}
  We study Bartnik's splitting conjecture through the space $\mathcal{C}_M$ of Cauchy hypersurfaces,
  equipped with the Hausdorff-type distance $d_J$ introduced by Lange and Peteranderl.
  For a timelike geodesically complete, spatially compact spacetime
  satisfying the timelike convergence condition, we show that the splitting of the spacetime
  is equivalent to each of the conditions that
  $(\mathcal{C}_M, d_J)$ contains a $(1, \varepsilon)$-quasi-geodesic line,
  is quasi-isometric to $\mathbf{R}$, is Gromov hyperbolic,
  and has at least two ends.
  Thus the Lorentzian problem of finding a timelike line becomes a problem in the coarse geometry of $(\mathcal{C}_M, d_J)$.
  The key step, which holds without the timelike convergence condition,
  is that any $(1, \varepsilon)$-quasi-geodesic line in $(\mathcal{C}_M, d_J)$ yields a complete timelike line in the spacetime.
  Moreover, we characterize the splitting by the topology of $\mathcal{C}_M$ alone,
  showing that the spacetime splits if and only if $\mathcal{C}_M$ is homeomorphic to $\mathbf{R}$ times the Hilbert cube.
  Hence Bartnik's conjecture is equivalent to the statement that this holds for every spacetime as above.
  We also show, under timelike geodesic completeness and spatial compactness alone,
  that $(\mathcal{C}_M, d_J)$ is a proper injective, hence geodesic, metric space.
  On the other hand, for every product spacetime with compact Riemannian factor, where it is Gromov hyperbolic,
  and for every spacetime satisfying the timelike convergence condition,
  it is neither CAT(0) nor Busemann convex.
\end{abstract}

\maketitle

\section{Introduction}\label{sec:intro}

A Lorentzian manifold does not carry a natural distance function,
and many tools of metric geometry
are not directly available in Lorentzian geometry.
Nevertheless, several metric spaces can be constructed from Lorentzian data,
and Lorentzian problems can sometimes be studied through them.
For instance, Sormani and Vega~\cite{SV} introduced the null distance
$\widehat{d}_\tau$ on a spacetime, built from a time function $\tau$,
and Burtscher and Garc\'ia-Heveling~\cite{BG} proved
that a spacetime is globally hyperbolic if and only if
it is metrically complete with respect to the null distance of a suitable time function,
which can be regarded as a Lorentzian substitute for the Hopf--Rinow theorem.
Similarly, Lange and Peteranderl~\cite{LP} introduced a Hausdorff-type metric $\dJ$
on the set $\CM$ of all Cauchy hypersurfaces of a globally hyperbolic spacetime $M$,
built from the Lorentzian distance (see~Definition~\ref{def:dJ}),
which is used in the quantitative Lorentzian isoperimetric inequality of~\cite[Theorem 1.6]{LP25}.

In this article, we study Bartnik's splitting conjecture, a rigidity problem in Lorentzian geometry,
from the viewpoint of the metric space $(\CM, \dJ)$
and of its underlying topological space.
Let us first recall the background.
In Riemannian geometry, a fundamental rigidity result is the splitting theorem of Cheeger and Gromoll.

\begin{fact}[Cheeger--Gromoll~\cite{CG}]
  Let $(N, h)$ be a complete Riemannian manifold with non-negative Ricci curvature.
  If $N$ contains a line, i.e., a geodesic $\gamma \colon \R \to N$ satisfying
  $\operatorname{dist}_h(\gamma(s), \gamma(t)) = \abs{t - s}$ for all $s, t \in \R$,
  then $(N, h)$ is isometric to a product $(\R \times N', dt^2 + h')$,
  where $(N', h')$ is a complete Riemannian manifold with non-negative Ricci curvature.
\end{fact}

The Lorentzian analogue of this theorem was conjectured by Yau~\cite[Problem 115]{Yau82}.
A spacetime satisfies the \emph{timelike convergence condition} if $\Ric(v, v) \geq 0$ for every timelike vector $v$,
and a \emph{complete timelike line} is a timelike geodesic $\gamma \colon \R \to M$ parametrized by proper time
with $d(\gamma(s), \gamma(t)) = t - s$ for all $s \leq t$, where $d$ is the Lorentzian distance.

\begin{fact}[Lorentzian splitting theorem~\cite{Es, Ga, Ne}]\label{fact:splitting}
  Let $(M, g)$ be a spacetime satisfying the timelike convergence condition
  which is globally hyperbolic or timelike geodesically complete.
  If $M$ contains a complete timelike line, then $(M, g)$ is isometric to
  $(\R \times S, -dt^2 + h)$ for a complete Riemannian manifold $(S, h)$ with non-negative Ricci curvature.
\end{fact}

Fact~\ref{fact:splitting} was proved by Eschenburg~\cite{Es} for globally hyperbolic, timelike geodesically complete spacetimes,
by Galloway~\cite{Ga} for globally hyperbolic spacetimes,
and by Newman~\cite{Ne} for timelike geodesically complete spacetimes.
In cosmology, one is interested in spatially compact spacetimes, i.e., spacetimes with compact Cauchy hypersurfaces.
For such spacetimes, Bartnik~\cite{Ba} proposed the following conjecture.

\begin{conjecture}[Bartnik's splitting conjecture~\cite{Ba}]\label{conj:bartnik}
  Let $(M, g)$ be a spatially compact, globally hyperbolic spacetime.
  If $(M, g)$ is timelike geodesically complete and satisfies the timelike convergence condition,
  then $(M, g)$ splits isometrically as $(\R \times S, -dt^2 + h)$,
  where $(S, h)$ is a compact Riemannian manifold with non-negative Ricci curvature.
\end{conjecture}

The conjecture is still open; see~\cite[Section 1.4]{Mi} for a survey.
By Fact~\ref{fact:splitting}, it suffices to find a complete timelike line in $M$.
This is a genuinely Lorentzian problem:
limits of maximizing curves may become null,
so that a line cannot simply be obtained as a limit of maximizing segments,
whereas in Riemannian geometry a complete manifold with at least two ends always contains a line obtained in this way.
In fact, without a curvature condition, a timelike line need not exist:
a spatially compact, globally hyperbolic, causally geodesically complete spacetime
which contains no timelike line was constructed in~\cite[Theorem 2]{EG};
since it does not satisfy the timelike convergence condition,
any construction of a timelike line has to use this condition.
Many partial results establish the splitting under additional conditions~\cite{Ba, Ga84, EG, GV, GV17, Ga19, Le19, CFH20, GaHe25, FHSo25};
we refer to the survey~\cite[Section 1.4]{Mi} and mention only some recent ones.
In three spacetime dimensions, the splitting was obtained under null geodesic completeness and the dominant energy condition,
except when the Cauchy hypersurfaces are spheres or projective planes~\cite{GLO26}.
Like the present article, several recent works reformulate the conjecture
in terms of Cauchy Busemann functions associated with a compact Cauchy hypersurface~\cite{Mi} (see~Fact~\ref{fact:defect}),
of Cauchy temporal functions for which the Lorentzian distances to and from the level sets are locally equi-Lipschitz~\cite{GMO},
and of $p$-harmonic temporal functions~\cite{Oh26}.
All these approaches work with objects in $M$, such as rays, horospheres, Busemann functions and temporal functions.
We propose a different viewpoint and transfer the problem to the space $\CM$ of all Cauchy hypersurfaces in two steps.

First, the Lorentzian problem of finding a timelike line becomes a problem in the coarse geometry
of the metric space $(\CM, \dJ)$ (Theorem~\ref{thm:main}).
Second, this problem is reformulated as a condition on the topology of $\CM$ alone (Corollary~\ref{cor:type}).

\emph{From Lorentzian geometry to coarse geometry.}
Our main result characterizes the splitting by the coarse geometry of $(\CM, \dJ)$.
The relevant notions of metric geometry are recalled in Section~\ref{subsec:metric}.

\begin{reptheorem}{Theorem}{thm:main}
  Let $(M, g)$ be a spatially compact, globally hyperbolic spacetime.
  Suppose that $(M, g)$ is timelike geodesically complete and satisfies
  the timelike convergence condition.
  Then the following conditions are equivalent:
  \begin{enumerate}
    \item $(M, g)$ splits isometrically as $(\R \times S, -dt^2 + h)$,
      where $(S, h)$ is a compact Riemannian manifold
      with non-negative Ricci curvature;
    \item $(\CM, \dJ)$ contains a $(1, \varepsilon)$-quasi-geodesic line
      for some $\varepsilon \geq 0$;
    \item $(\CM, \dJ)$ is quasi-isometric to $\R$;
    \item $(\CM, \dJ)$ is Gromov hyperbolic;
    \item $(\CM, \dJ)$ has at least two ends.
  \end{enumerate}
\end{reptheorem}

\begin{corollary}\label{cor:coarse}
  Bartnik's conjecture holds if and only if,
  for every spatially compact, globally hyperbolic spacetime $(M, g)$
  which is timelike geodesically complete and satisfies the timelike convergence condition,
  the space $(\CM, \dJ)$ is Gromov hyperbolic,
  or equivalently, quasi-isometric to $\R$,
  or equivalently, has at least two ends.
\end{corollary}

Thus Bartnik's conjecture becomes a problem in the coarse geometry
of a single metric space attached to the spacetime,
and the tools of coarse geometry, such as quasi-isometry invariance, become available.
The Lorentzian difficulty described above, namely that limits of maximizing curves may become null,
is handled once and for all by Theorem~\ref{thm:quasi} below,
which turns quasi-geodesic lines in $(\CM, \dJ)$ into timelike lines in $M$.

\emph{From coarse geometry to topology.}
Condition (v) is in fact topological, since the ends of a space depend only on its topology.
By Lemma~\ref{lem:topology}, the topology of $\dJ$ is the topology of uniform convergence of the graph functions of Cauchy hypersurfaces
with respect to a diffeomorphism $M \cong \R \times S'$ given by a Cauchy temporal function.
This is the $C^0$ topology of graphs, which does not involve the Lorentzian distance.
The following corollary characterizes the splitting by the topology of $\CM$ alone.
Let $Q \coloneqq [0, 1]^{\N}$ be the Hilbert cube, i.e., the product of countably many copies of $[0, 1]$ with the product topology,
and for $p \in M$, let $\CM(p)$ be the set of Cauchy hypersurfaces passing through $p$.

\begin{reptheorem}{Corollary}{cor:type}
  Under the assumptions of Theorem~\ref{thm:main}, the following conditions are equivalent:
  \begin{enumerate}
    \item $(M, g)$ splits isometrically as $(\R \times S, -dt^2 + h)$;
    \item $\CM$ is homeomorphic to $\R \times Q$;
    \item $\CM(p)$ is compact for some $p \in M$;
    \item $\CM(p)$ is homeomorphic to $Q$ for every $p \in M$.
  \end{enumerate}
  Moreover, if $(M, g) = (\R \times S, -dt^2 + h)$ and $p = (t_0, x_0)$,
  then there exists a homeomorphism of pairs $(\CM, \CM(p)) \cong (\R \times Q, \{t_0\} \times Q)$.
\end{reptheorem}
Consequently, Bartnik's conjecture is equivalent to the statement that $\CM$ is homeomorphic to $\R \times Q$
for every spacetime as in Theorem~\ref{thm:main};
in the split case, $\CM(p)$ corresponds to a slice $\{t_0\} \times Q$ separating the two ends.
Condition (iii) is equivalent to the point condition of Bartnik (Remark~\ref{rem:point}),
under which he obtained the splitting through a constant mean curvature Cauchy hypersurface~\cite[Theorem and Corollary 2]{Ba}.
Thus (iii) $\Rightarrow$ (i) recovers his result,
and the new content of Corollary~\ref{cor:type} is the determination of the homeomorphism types
of $\CM$ and of the pair $(\CM, \CM(p))$.

By Proposition~\ref{prop:two-ends}, without the timelike convergence condition or timelike geodesic completeness,
$\CM$ has one or two ends,
and it has two ends if and only if no Cauchy hypersurface crosses a sufficiently thick time band.
Without the timelike convergence condition, this condition lies strictly between the point condition and the min-min condition of~\cite{GaHe25}
(Proposition~\ref{prop:tip}, Remark~\ref{rem:causal} and Example~\ref{ex:cylinder});
in particular, (v) $\Rightarrow$ (i) of Theorem~\ref{thm:main} can also be derived from~\cite[Theorem 1.4]{GaHe25}.

The key step of the proof of Theorem~\ref{thm:main} is the following result,
which does not use any curvature condition.

\begin{reptheorem}{Theorem}{thm:quasi}
  Let $(M, g)$ be a timelike geodesically complete, spatially compact, globally hyperbolic spacetime,
  and let $(S_t)_{t \in \R}$ be a $(1, \varepsilon)$-quasi-geodesic line in $(\CM, \dJ)$.
  Then there exists a complete timelike line $\gamma$ passing through a point $z \in S_0$.
  Moreover, if $\gamma$ is parametrized by proper time with $\gamma(0) = z$, then
  \begin{equation}
    d(S_0, \gamma(t)) \leq t + 3\varepsilon \quad (t \geq 0), \qquad
    d(\gamma(t), S_0) \leq \abs{t} + 3\varepsilon \quad (t \leq 0).
  \end{equation}
\end{reptheorem}

We make three remarks on Theorem~\ref{thm:quasi}.
First, the additive control of the quasi-geodesic line gives $d(S_0, \gamma(t)) \leq d(z, \gamma(t)) + 3\varepsilon$ for $t > 0$
along the limit $\gamma$ of maximizing segments,
which ties the length of $\gamma$ to the Lorentzian distance from $S_0$ and thereby rules out a null limit.
The estimate in Theorem~\ref{thm:quasi} then follows by parametrizing the timelike line $\gamma$ by proper time.
Second, the quasi-geodesic line is arbitrary;
in particular, it is not assumed to be causally monotone in the sense of Definition~\ref{def:monotone}.
This is essential, since quasi-geodesic lines obtained from Gromov hyperbolicity, quasi-isometries or ends
carry no information on causality.
Third, the multiplicative constant $1$ cannot be replaced by any $\lambda \geq \sqrt{2}$
(we do not know whether it can be replaced by some $\lambda \in (1, \sqrt{2})$),
since suitably rescaled distance curves of Cauchy hypersurfaces are always $(\sqrt{2}, 0)$-quasi-geodesic lines,
even in the example of~\cite[Theorem 2]{EG}, which contains no timelike line (Proposition~\ref{prop:lambda}).
For this example, $(\CM, \dJ)$ is not Gromov hyperbolic and has exactly one end (Example~\ref{ex:EG}).

Theorem~\ref{thm:main} fails without the timelike convergence condition.
For warped products $-dt^2 + a(t)^2 h$ with $a$ bounded above and below by positive constants,
$(\CM, \dJ)$ is quasi-isometric to $\R$ and Gromov hyperbolic, although the spacetime need not split (Example~\ref{ex:warped}).
For the de Sitter spacetime $\mathbf{S}^{n+1}_1$, the space $(\CM, \dJ)$ contains
$\mathbf{H}^{n+1} \times \R$ with the $\ell^1$-product distance isometrically,
so that it contains geodesic lines but is not Gromov hyperbolic (Example~\ref{ex:desitter}).

The proof of the characterization by ends uses that $(\CM, \dJ)$ is a geodesic space.
Lange and Peteranderl proved the geodesic property under timelike Cauchy completeness~\cite{LP}.
We prove it without this assumption, using only the standing assumption of timelike geodesic completeness
through the compactness of the parallel sets of a compact Cauchy hypersurface~\cite[Theorem 2.4]{Mi}.
The key observation is that closed balls are order intervals with respect to the causal order on $\CM$ (Proposition~\ref{prop:ball}),
which yields the following Helly-type property.

\begin{reptheorem}{Theorem}{thm:injective}
  Let $(S_i, r_i)_{i \in I}$ be a family in $\CM \times \R$ with $\dJ(S_i, S_j) \leq r_i + r_j$ for all $i, j \in I$.
  Then $\bigcap_{i \in I} \overline{B}(S_i, r_i) \neq \emptyset$.
  Consequently, $(\CM, \dJ)$ is a hyperconvex, and hence injective, metric space.
  In particular, it is a proper geodesic metric space.
\end{reptheorem}

In the proof of Theorem~\ref{thm:main}, we only use the consequence that $(\CM, \dJ)$ is a proper geodesic metric space;
the injectivity is an independent structural result.

Finally, injective metric spaces share several features with non-positively curved spaces~\cite{La}.
However, $(\CM, \dJ)$ is neither a CAT(0) space nor Busemann convex
as soon as some smooth, spacelike, acausal compact Cauchy hypersurface $\Sigma$ satisfies $D_\infty(\Sigma) > 0$ (Theorem~\ref{thm:convexity}),
in particular for every product spacetime and, under the timelike convergence condition, for every spacetime (Corollary~\ref{cor:convexity}),
whereas the $L^2$-metric of~\cite{Mon} has non-positive sectional curvature.

The article is organized as follows.
In Section~\ref{sec:prelim}, we fix our conventions on spacetimes,
recall the metric $\dJ$ and some notions of metric geometry,
and collect the results of~\cite{Mi} which will be used.
In Section~\ref{sec:basic}, we prove basic properties of $\dJ$,
including the properness of $(\CM, \dJ)$,
show that the topology of $\dJ$ is the $C^0$ topology of graphs,
and introduce the signed distance and the causal order on $\CM$.
In Section~\ref{sec:injective}, we prove Theorem~\ref{thm:injective}.
In Section~\ref{sec:parallel}, we study distance curves.
In Section~\ref{sec:line}, we prove Theorem~\ref{thm:quasi}.
In Section~\ref{sec:main}, we study product spacetimes and Gromov hyperbolicity,
prove Theorem~\ref{thm:main}
and Corollary~\ref{cor:type},
and give the examples without the timelike convergence condition.
In Section~\ref{sec:convexity}, we prove Theorem~\ref{thm:convexity}
and compare $\dJ$ with the $L^2$-metric of~\cite{Mon}.
In Section~\ref{sec:remarks}, we give
concluding remarks.

\section{Preliminaries}\label{sec:prelim}

In this section, we fix our conventions and recall the results which will be used
in the proofs of the main theorems.

\subsection{Spacetimes}

A \emph{spacetime} $(M, g)$ is a smooth, connected, time-oriented Lorentzian manifold of dimension $n + 1$, $n \geq 1$,
with signature $(-, +, \dots, +)$.
For causality theory, we follow the conventions of the review~\cite{MiLR},
which are also those of~\cite[Section 1.5]{Mi}, with two differences.
First, causal curves may be past-directed, whereas causal curves in~\cite{MiLR, Mi} are future-directed;
statements of~\cite{MiLR, Mi} about causal curves refer to future-directed ones.
Second, Cauchy hypersurfaces (Definition~\ref{def:cauchy}) are not assumed to be acausal.

A \emph{future-directed} (resp.\ \emph{past-directed}) \emph{causal curve} is a piecewise $C^1$ curve $\gamma \colon I \to M$,
defined on an interval $I \subset \R$ containing at least two points,
whose tangent vectors, including both one-sided tangent vectors at break points,
are future-directed (resp.\ past-directed) causal (see~\cite[Section 1.11]{MiLR}).
Future-directed and past-directed \emph{timelike curves} are defined analogously,
and a \emph{causal} (resp.\ \emph{timelike}) \emph{curve} is a future-directed or past-directed one.
For $p, q \in M$, we write $p \ll q$ (resp.\ $p \leq q$) if there exists a future-directed timelike
(resp.\ causal) curve from $p$ to $q$, or $p = q$ in the case of $\leq$, and we write $p < q$ if $p \leq q$ and $p \neq q$.
For $A \subset M$, we set $I^{+}(A) \coloneqq \{q \in M ; p \ll q \text{ for some } p \in A\}$ and
$J^{+}(A) \coloneqq \{q \in M ; p \leq q \text{ for some } p \in A\}$, and define $I^-(A)$ and $J^-(A)$ similarly.
A future-directed causal curve $\gamma \colon I \to M$ has a \emph{future endpoint} $p \in M$
if $\gamma(t) \to p$ as $t \to \sup I$ with $t < \sup I$,
and it is \emph{future inextendible} if it has no future endpoint (in particular, $\sup I \notin I$);
past endpoints and past inextendibility are defined with $\inf I$,
and the roles of $\sup I$ and $\inf I$ are exchanged for past-directed curves.
A causal curve is \emph{inextendible} if it is both future and past inextendible.
The limit curve theorems~\cite[Section 2.11]{MiLR} produce \emph{continuous causal curves} in the sense of~\cite[Section 2.5]{MiLR},
for which these notions are defined in the same way~\cite[Definition 2.16]{MiLR}.
Using continuous causal curves instead of causal curves does not change the relation $\leq$~\cite[Remark 2.14]{MiLR}.

\begin{definition}[{\cite[Section 2.9]{MiLR}, \cite[Section 1.5]{Mi}}]\label{def:distance}
  The \emph{Lorentzian distance} $d \colon M \times M \to [0, \infty]$ is defined by
  \begin{equation}
    d(p, q) \coloneqq
    \begin{cases}
      \displaystyle \sup_{\gamma} L_g(\gamma), & q \in J^+(p), \\
      0,                                       & q \notin J^+(p),
    \end{cases}
  \end{equation}
  where the supremum is taken over all future-directed continuous causal curves $\gamma$ from $p$ to $q$,
  and $L_g(\gamma)$ denotes the Lorentzian length of $\gamma$~\cite[Section 2.10]{MiLR}.
  For subsets $A, B \subset M$, we set
  \begin{equation}
    d(A, B) \coloneqq \sup_{p \in A, \, q \in B} d(p, q).
  \end{equation}
\end{definition}

The Lorentzian distance does not depend on the regularity of the curves.
Indeed, the supremum is unchanged if it is taken over causal curves only,
since every continuous causal curve from $p$ to $q$ admits an interpolating piecewise causal geodesic
from $p$ to $q$ which is at least as long~\cite[Theorem 2.37]{MiLR}.
We note that some authors, especially in synthetic Lorentzian geometry,
set $d(p, q) = -\infty$ if $q \notin J^+(p)$; we do not adopt this convention, so that $d \geq 0$.
The Lorentzian distance satisfies the \emph{reverse triangle inequality}
$d(p, r) \geq d(p, q) + d(q, r)$ for all $p \leq q \leq r$~\cite[Theorem 2.32]{MiLR},
and $d(p, q) > 0$ if and only if $p \ll q$.
A \emph{time function} is a continuous function $t \colon M \to \R$ with $t(p) < t(q)$ for all $p < q$,
and it is a \emph{Cauchy time function} if each of its level sets is met by every inextendible causal curve~\cite[Definitions 1.21 and 1.29]{MiLR}.
A \emph{Cauchy temporal function} is a smooth Cauchy time function with past-directed timelike gradient.

\begin{definition}[{\cite[Definition 14.28]{On}, cf.~\cite[Definition 3.35]{MiLR}}]\label{def:cauchy}
  A \emph{Cauchy hypersurface} of $M$ is a subset $S \subset M$ which is met exactly once
  by every inextendible timelike curve.
  A Cauchy hypersurface $S$ is called \emph{acausal} if no causal curve meets $S$ twice.
  A spacetime is called \emph{globally hyperbolic} if it admits a Cauchy hypersurface,
  and \emph{spatially compact} if it admits a compact Cauchy hypersurface.
  We denote by $\CM$ the set of all Cauchy hypersurfaces of $M$.
\end{definition}

Definition~\ref{def:cauchy} does not change if locally Lipschitz or continuous timelike curves are used,
since $\ll$ is the same for piecewise $C^1$ and locally Lipschitz timelike curves~\cite[Theorem 2.15 and the discussion after Corollary 2.16]{GKSS}.
Let $S \in \CM$ and let $\gamma$ be an inextendible continuous timelike curve.
If $\gamma$ met $S$ twice, joining the two points by a timelike curve and extending it by integral curves of a timelike vector field
would give an inextendible timelike curve meeting $S$ twice.
If $\gamma$ missed $S$, then $\gamma \subset I^+(S)$, say, and joining points $\gamma(t_n)$, $t_n \searrow \inf I$, by timelike curves
gives a timelike curve missing $S$, which is past inextendible
since $M$ is globally hyperbolic (see below) and hence $\gamma$ is not partially imprisoned~\cite[Theorem 2.80]{MiLR}.

Every $S \in \CM$ is \emph{achronal}, i.e., there are no $p, q \in S$ with $p \ll q$, it is closed, and it satisfies $D(S) = M$,
where $D(S)$ is the Cauchy development defined by inextendible causal curves~\cite[Definition 3.29]{MiLR}
(Lemma~\ref{lem:cauchy} (iv) below, whose proof does not use global hyperbolicity).
Hence, if $M$ admits some $S \in \CM$, then $\operatorname{Int} D(S) = M$,
and $M$ is globally hyperbolic in the sense of~\cite{MiLR} by~\cite[Theorem 3.45]{MiLR}.
Conversely, a spacetime which is globally hyperbolic in the sense of~\cite{MiLR}
admits a Cauchy hypersurface in the sense of~\cite[Definition 3.35]{MiLR} by~\cite[Theorem 4.119]{MiLR},
which clearly belongs to $\CM$.
Thus the two notions of global hyperbolicity agree.
By the theorem of Bernal and S\'anchez~\cite{BS}, \cite[Theorem 4.120]{MiLR},
$M$ admits a Cauchy temporal function $\mathcal{T}$,
whose level sets are smooth spacelike, hence acausal, elements of $\CM$.
All elements of $\CM$ are homeomorphic to each other by Geroch~\cite[Property 7]{Ge}.
Hence if $M$ is spatially compact, then every element of $\CM$ is compact,
and in particular the level sets of $\mathcal{T}$ are compact acausal Cauchy hypersurfaces.
Finally, the acausal elements of $\CM$ are exactly the Cauchy hypersurfaces in the sense of~\cite[Definition 3.35]{MiLR},
i.e., the closed acausal sets $S$ with $D(S) = M$, which are those considered in~\cite{Mi}.
Indeed, an acausal $S \in \CM$ is such a set by Lemma~\ref{lem:cauchy} (iv),
and conversely, such a set is met by every inextendible timelike curve since $D(S) = M$, and at most once by acausality.

\begin{fact}[{\cite[Theorems 4.123 and 4.124]{MiLR}}]\label{fact:gh}
  Let $(M, g)$ be a globally hyperbolic spacetime.
  Then $d$ is finite and continuous on $M \times M$.
  Moreover, for all $p, q \in M$ with $p < q$, there exists a maximizing causal geodesic
  from $p$ to $q$, i.e., a causal geodesic of length $d(p, q)$ (Avez~\cite{Av} and Seifert~\cite{Se}).
\end{fact}

\begin{definition}[{\cite[Section 1.5]{Mi}}]\label{def:line}
  A \emph{timelike line} is an inextendible timelike geodesic $\gamma$ such that
  every segment of $\gamma$ is maximizing for $d$.
  A timelike line is called \emph{complete} if its proper time parameter ranges over all of $\R$.
\end{definition}

A spacetime is called \emph{timelike geodesically complete} if every inextendible timelike geodesic
is defined on all of $\R$ when parametrized by proper time.
The timelike convergence condition and the Lorentzian splitting theorem (Fact~\ref{fact:splitting}) were recalled in Section~\ref{sec:intro}.

We call $(\R \times S, -dt^2 + h)$, where $(S, h)$ is a Riemannian manifold, a \emph{product spacetime}.
It satisfies the timelike convergence condition if and only if $(S, h)$ has non-negative Ricci curvature,
since $\Ric(\partial_t, \cdot) = 0$ and $\Ric(w, w) = \Ric_h(w, w)$ for $w$ tangent to $S$.
If $M$ is spatially compact and splits as in Fact~\ref{fact:splitting},
then $S$ is compact, since $\{0\} \times S$ is a Cauchy hypersurface.
In this case, we say that $(M, g)$ \emph{splits}.

\emph{Throughout the paper, $(M, g)$ denotes a globally hyperbolic, spatially compact spacetime
which is timelike geodesically complete.
The timelike convergence condition is assumed only when explicitly stated.}

These assumptions and the timelike convergence condition are invariant under reversing the time orientation of $(M, g)$.
This reversal interchanges future-directed and past-directed curves, $I^+$ and $I^-$, and $J^+$ and $J^-$,
replaces $d(p, q)$ by $d(q, p)$, and does not change the Cauchy hypersurfaces.
For a statement about $(M, g)$, its \emph{time-dual statement} is the statement obtained by applying it
to $(M, g)$ with the reversed time orientation and rewriting the result in terms of $(M, g)$.
Hence the time-dual statement of every result proved under our assumptions also holds.

\subsection{The space of Cauchy hypersurfaces}

The following metric on $\CM$ was introduced in~\cite{LP},
based on an idea of Bahn and Ehrlich~\cite{BE99}.

\begin{definition}[{\cite{LP}}]\label{def:dJ}
  For $A, B \subset M$, we set
  \begin{equation}
    \dJ(A, B) \coloneqq \sup_{x \in A, \, y \in B} \bigl( d(x, y) + d(y, x) \bigr).
  \end{equation}
\end{definition}

The article~\cite{LP} uses locally Lipschitz curves.
Since causal curves are locally Lipschitz, and locally Lipschitz causal curves are continuous causal curves with the same length~\cite[Theorem 2.12 and Section 2.10]{MiLR},
the remarks after Definitions~\ref{def:distance} and~\ref{def:cauchy} show that $d$, $\dJ$ and $\CM$ are the same in both conventions.
Thus the results of~\cite{LP} are unaffected by our convention.

\begin{fact}[{\cite[Theorem 1.1 (iii), (iv)]{LP}}]\label{fact:LP}
  Let $(M, g)$ be a globally hyperbolic spacetime.
  Then $\dJ$ is an extended metric on $\CM$, i.e., a metric which may take the value $\infty$.
  Moreover, the following hold:
  \begin{enumerate}
    \item if $(M, g)$ is spatially compact, then $(\CM, \dJ)$ is a locally compact metric space;
    \item if $(M, g)$ is spatially compact and timelike Cauchy complete,
          then $(\CM, \dJ)$ is a proper geodesic metric space.
  \end{enumerate}
\end{fact}

Since $(M, g)$ is spatially compact throughout the paper, $(\CM, \dJ)$ is a metric space.
Timelike Cauchy completeness, introduced by Busemann~\cite{Bu67} (see also~\cite{Be76, Ta}), will not be used:
$(\CM, \dJ)$ is proper and injective, and hence geodesic, under our standing assumptions alone
(Proposition~\ref{prop:proper} and Theorem~\ref{thm:injective}).

\subsection{Notions from metric geometry}\label{subsec:metric}

Let $(X, \rho)$ be a metric space.

\begin{definition}[{\cite[Chapter I.8, Definitions 8.14 and 8.22]{BH}}]\label{def:quasi}
  Let $\lambda \geq 1$ and $\varepsilon \geq 0$.
  A family $(x_t)_{t \in \R}$ of points in $X$ is called a \emph{$(\lambda, \varepsilon)$-quasi-geodesic line} if
  \begin{equation}
    \frac{1}{\lambda} \abs{t - s} - \varepsilon \leq \rho(x_s, x_t) \leq \lambda \abs{t - s} + \varepsilon \qquad \text{for all } s, t \in \R.
  \end{equation}
  No continuity in $t$ is required.
  A $(1, 0)$-quasi-geodesic line is called a \emph{geodesic line}.
\end{definition}

In the terminology of~\cite{BH}, a $(\lambda, \varepsilon)$-quasi-geodesic line is a $(\lambda, \varepsilon)$-quasi-isometric embedding of $\R$.
We note that a geodesic line in a metric space is different from a timelike line in a spacetime.

\begin{definition}[{\cite[Chapter I.8, Definition 8.14]{BH}}]\label{def:rough}
  Let $(Y, \rho')$ be a metric space, $\lambda \geq 1$ and $\varepsilon \geq 0$.
  A map $\varphi \colon X \to Y$ is called a \emph{$(\lambda, \varepsilon)$-quasi-isometry} if
  \begin{equation}
    \frac{1}{\lambda} \rho(x, x') - \varepsilon \leq \rho'(\varphi(x), \varphi(x')) \leq \lambda \rho(x, x') + \varepsilon
    \qquad \text{for all } x, x' \in X,
  \end{equation}
  and there exists $C \geq 0$ such that every point of $Y$ lies within distance $C$ of $\varphi(X)$.
  We say that $X$ is \emph{quasi-isometric} to $Y$ if there exists a $(\lambda, \varepsilon)$-quasi-isometry
  $\varphi \colon X \to Y$ for some $\lambda \geq 1$ and $\varepsilon \geq 0$.
\end{definition}

\begin{definition}[{\cite[Chapter III.H, Definitions 1.19 and 1.20]{BH}}]\label{def:hyperbolic}
  The \emph{Gromov product} of $y, z \in X$ with respect to $x \in X$ is
  \begin{equation}
    (y \cdot z)_x \coloneqq \frac{1}{2} \bigl( \rho(y, x) + \rho(z, x) - \rho(y, z) \bigr).
  \end{equation}
  For $\delta \geq 0$, the space $X$ is called \emph{$\delta$-hyperbolic} if
  \begin{equation}
    (x \cdot y)_w \geq \min\{(x \cdot z)_w, (z \cdot y)_w\} - \delta \qquad \text{for all } w, x, y, z \in X,
  \end{equation}
  and \emph{Gromov hyperbolic} if it is $\delta$-hyperbolic for some $\delta \geq 0$.
\end{definition}

We note that Definition~\ref{def:hyperbolic} does not require $X$ to be a geodesic space.

\begin{definition}[{\cite[Chapter I.8, Definition 8.27]{BH}}]\label{def:ends}
  A \emph{ray} in $X$ is a continuous map $r \colon [0, \infty) \to X$,
  and it is \emph{proper} if $r^{-1}(K)$ is compact for every compact set $K \subset X$.
  Two proper rays $r_1$ and $r_2$ \emph{converge to the same end} if for every compact set $K \subset X$
  there exists $N \geq 0$ such that $r_1([N, \infty))$ and $r_2([N, \infty))$ lie in the same path component of $X \setminus K$.
  This is an equivalence relation on proper rays, and its equivalence classes are called the \emph{ends} of $X$.
  We say that $X$ has \emph{at least two ends} (resp.\ \emph{exactly one end}) if the number of ends of $X$ is at least two
  (resp.\ exactly one).
\end{definition}

For instance, a geodesic ray $r \colon [0, \infty) \to X$, i.e., an isometric embedding, is a proper ray.

\begin{definition}[{\cite[Chapter I.1, Definition 1.3]{BH}, \cite[Definition 8.1.1]{Pa}}]\label{def:local-geodesic}
  A curve $c \colon I \to X$ defined on an interval $I \subset \R$ is called a \emph{local geodesic}
  if every $t_0 \in I$ has a neighborhood $J \subset I$ such that $\rho(c(s), c(t)) = \abs{t - s}$ for all $s, t \in J$.
  A geodesic space $X$ is called \emph{Busemann convex}, or a \emph{Busemann space}, if,
  for any two affinely parametrized geodesics $\gamma \colon [a, b] \to X$ and $\gamma' \colon [a', b'] \to X$,
  the function $[a, b] \times [a', b'] \to \R$, $(t, t') \mapsto \rho(\gamma(t), \gamma'(t'))$, is convex.
  Here a curve is an affinely parametrized geodesic if it is a geodesic parametrized proportionally to arc length.
\end{definition}

\begin{fact}[{\cite[Chapter II.1, Proposition 1.4 (2)]{BH}, \cite[Corollary 8.2.3]{Pa}}]\label{fact:cat0}
  In a CAT(0) space, and in a Busemann convex space, every local geodesic is a geodesic.
\end{fact}

\subsection{Cauchy Busemann functions}\label{subsec:busemann}

We recall the notation and some results of~\cite{Mi};
numbering refers to the second arXiv version (arXiv:2610.17092v2).
In~\cite{Mi}, Cauchy hypersurfaces are assumed to be acausal, ``timelike complete''
means timelike geodesically complete, and the background causality theory is that of~\cite{MiLR} (see~\cite[Section 1.5]{Mi}).
Let $\Sigma$ be an acausal compact Cauchy hypersurface.
We set $d^+_\Sigma(q) \coloneqq d(\Sigma, q)$ and $d^-_\Sigma(q) \coloneqq d(q, \Sigma)$ for $q \in M$,
and for $c > 0$, we define the \emph{parallel sets}~\cite[Section 5.5]{LP}
\begin{equation}
  \Sigma_c \coloneqq (d^+_\Sigma)^{-1}(c), \qquad \Sigma_{-c} \coloneqq (d^-_\Sigma)^{-1}(c),
\end{equation}
and $\Sigma_0 \coloneqq \Sigma$.
The \emph{signed distance} from $\Sigma$ is the function $d_\Sigma \coloneqq d^+_\Sigma - d^-_\Sigma \colon M \to \R$~\cite[Proposition 2.9]{Mi20}, \cite[Section 2]{Mi}.
Since $\Sigma$ is achronal, $d^-_\Sigma = 0$ on $J^+(\Sigma)$ and $d^+_\Sigma = 0$ on $J^-(\Sigma)$,
so that $d_\Sigma = d^+_\Sigma$ on $J^+(\Sigma)$ and $d_\Sigma = -d^-_\Sigma$ on $J^-(\Sigma)$.
For $\tau > 0$, the \emph{Cauchy Busemann functions} are defined by
\begin{equation}
  \begin{gathered}
    b^+_\tau \coloneqq \tau - d(\cdot, \Sigma_\tau), \qquad b^+ \coloneqq \lim_{\tau \to \infty} b^+_\tau, \\
    b^-_\tau \coloneqq -\tau + d(\Sigma_{-\tau}, \cdot), \qquad b^- \coloneqq \lim_{\tau \to \infty} b^-_\tau.
  \end{gathered}
\end{equation}
We write $\dom b^+ \coloneqq \{b^+ < \infty\}$ and $\dom b^- \coloneqq \{b^- > -\infty\}$.
When the dependence on $\Sigma$ matters, we write $b^\pm_\Sigma$.
Under the reversal of the time orientation, $d^\pm_\Sigma$, $\Sigma_c$, $d_\Sigma$ and $b^\pm$ are replaced by
$d^\mp_\Sigma$, $\Sigma_{-c}$, $-d_\Sigma$ and $-b^\mp$, respectively, while $\CM$ and $\dJ$ are unchanged.

\begin{fact}[{\cite[Theorem 2.1]{Mi}}]\label{fact:cont}
  The function $d^+_\Sigma$ is finite and continuous on $M$, and so is $d^-_\Sigma$.
\end{fact}

\begin{fact}[{\cite[Theorem 2.4]{Mi}}]\label{fact:levels}
  For each $c \in \R$, the set $\Sigma_c$ is an acausal compact Cauchy hypersurface,
  and for $c \geq 0$ the set $J^+(\Sigma) \cap J^-(\Sigma_c)$ is compact.
  The function $d_\Sigma$ is a Cauchy time function whose level sets are $\Sigma_c$, and
  $(d^+_{\Sigma_a})^{-1}(b) = \Sigma_{a + b}$ for all $a, b > 0$.
  In particular, $d^+_\Sigma = d_\Sigma \to \infty$ along every future inextendible causal curve,
  since $d_\Sigma$ is a Cauchy time function.
  The time-dual statements also hold.
\end{fact}

\begin{fact}[{\cite[Proposition 2.6 and Corollary 2.7]{Mi}}]\label{fact:lengths}
  For $a \leq b$, every causal curve from $\Sigma_a$ to $\Sigma_b$ has length at most $b - a$;
  in particular, $d(\Sigma_a, \Sigma_b) \leq b - a$.
  For $0 < \tau \leq \tau'$ and $y \in M$, we have
  \begin{equation}
    d(y, \Sigma_\tau) \leq d(y, \Sigma_{\tau'}) \leq d(y, \Sigma_\tau) + \tau' - \tau.
  \end{equation}
  The time-dual statements also hold.
\end{fact}

\begin{fact}[{\cite[Theorems 7.2 and 7.3]{Mi}}]\label{fact:defect}
  The sequence $\Delta_k \coloneqq 2k - d(\Sigma_{-k}, \Sigma_k)$, $k \in \N$, is non-negative and non-decreasing, and
  \begin{equation}
    \lim_{k \to \infty} \Delta_k = \min_\Sigma (b^+ - b^-) = \inf_M (b^+ - b^-) \in [0, \infty].
  \end{equation}
  The following conditions are equivalent:
  \begin{enumerate}
    \item $\sup_k \Delta_k < \infty$;
    \item $\dom b^+ \cap \dom b^- \neq \emptyset$;
    \item $b^+ - b^-$ attains a finite minimum.
  \end{enumerate}
  Under the timelike convergence condition, these conditions are equivalent to the splitting of $(M, g)$.
\end{fact}

\begin{fact}[{\cite[Example 7.4]{Mi}}]\label{fact:product}
  If $(M, g) = (\R \times S, -dt^2 + h)$, where $(S, h)$ is a compact connected Riemannian manifold, and $\Sigma = \{t = f(x)\}$,
  then $\lim_{k \to \infty} \Delta_k = \max f - \min f$.
\end{fact}

\section{Basic properties of the metric $\dJ$}\label{sec:basic}

In this section, we prove some elementary properties of $\dJ$ which will be used repeatedly,
show that $(\CM, \dJ)$ is proper under our standing assumptions,
and introduce the signed distance and the causal order on $\CM$.

\begin{lemma}\label{lem:max}
  For all $x, y \in M$, at least one of $d(x, y)$ and $d(y, x)$ vanishes.
  Consequently, $\dJ(A, B) = \max\{d(A, B), d(B, A)\}$ for all $A, B \subset M$.
\end{lemma}

\begin{proof}
  If both $d(x, y)$ and $d(y, x)$ were positive, we would have $x \ll y \ll x$.
  This contradicts the chronology condition, which holds since $M$ is globally hyperbolic.
\end{proof}

\begin{lemma}[{cf.~\cite[Lemma 14.29]{On}}]\label{lem:cauchy}
  Let $S \in \CM$. Then the following hold:
  \begin{enumerate}
    \item $M$ is the disjoint union of $I^-(S)$, $S$ and $I^+(S)$.
          In particular, $J^+(S) \setminus S = I^+(S)$ and $J^+(S) \cap I^-(S) = \emptyset$;
    \item every causal curve from a point of $I^-(S)$ to a point of $J^+(S)$ meets $S$;
    \item the function $q \mapsto d(S, q)$ is continuous on $M$, and the supremum defining it is attained;
    \item $S$ is closed, and every inextendible causal curve meets $S$.
  \end{enumerate}
  The time-dual statements of (ii) and (iii) also hold.
\end{lemma}

\begin{proof}
  (i) Let $p \in M$. An inextendible timelike curve through $p$ meets $S$ at exactly one point $s$,
  and either $s = p$, $s \ll p$ or $p \ll s$.
  Thus $M = I^-(S) \cup S \cup I^+(S)$.
  The three sets are disjoint by the achronality of $S$.
  Since $S$ is achronal and $M = I^-(S) \cup S \cup I^+(S)$,
  the set $S$ is closed and maximal achronal, and $J^\pm(S) \setminus I^\pm(S) = S$~\cite[Theorem 2.146]{MiLR},
  i.e., $J^\pm(S) = S \cup I^\pm(S)$.
  This gives the remaining statements.

  (ii) By (i), $M \setminus S$ is the disjoint union of the open sets $I^-(S)$ and $I^+(S)$, and $J^+(S) \setminus S = I^+(S)$.
  Hence a causal curve from a point of $I^-(S)$ to a point of $J^+(S)$ which does not meet $S$
  would have a connected image covered by two disjoint open sets meeting it, which is impossible.

  (iii) This follows from the continuity of $d$ (Fact~\ref{fact:gh}) and the compactness of $S$.

  (iv) The set $S$ is closed by the proof of (i).
  We recall from~\cite[Definition 3.29]{MiLR} the Cauchy developments
  \begin{equation}
    \begin{aligned}
      D(S) & = \{p \in M ; \text{every inextendible causal curve through } p \text{ meets } S\}, \\
      \tilde{D}(S) & = \{p \in M ; \text{every inextendible timelike curve through } p \text{ meets } S\},
    \end{aligned}
  \end{equation}
  where, following~\cite{MiLR}, the curves are future-directed and piecewise $C^1$.
  This is no restriction, since reversing the parametrization of an inextendible past-directed causal (resp.\ timelike) curve
  gives an inextendible future-directed one with the same image.
  Thus $D(S) = M$ means that every inextendible causal curve in our sense meets $S$,
  and $\tilde{D}(S) = M$ means the same for inextendible timelike curves.
  By Definition~\ref{def:cauchy}, we have $\tilde{D}(S) = M$.
  Since $S$ is closed and achronal, \cite[Corollary 3.34]{MiLR} gives $D(S) = M$,
  i.e., every inextendible causal curve meets $S$.
\end{proof}

For acausal $S, S'$, the following lemma is~\cite[Lemma 3.8]{Mi},
since the constant $R(S, S')$ there coincides with $\dJ(S, S')$ by Lemma~\ref{lem:max}.
The proof in~\cite{Mi} uses acausality to see that $\partial[J^+(S) \cup J^+(S')]$ is a Cauchy hypersurface in the sense of~\cite{Mi};
we include a direct proof.

\begin{lemma}\label{lem:lipschitz}
  For all $S, S' \in \CM$ and $p \in M$, we have
  \begin{equation}
    \abs{d(S, p) - d(S', p)} \leq \dJ(S, S'), \qquad \abs{d(p, S) - d(p, S')} \leq \dJ(S, S').
  \end{equation}
\end{lemma}

\begin{proof}
  Let $R \coloneqq \dJ(S, S')$. We show $d(S, p) \leq d(S', p) + R$;
  the reverse inequality follows by exchanging $S$ and $S'$,
  and the second inequality is the time-dual statement.
  We may assume $d(S, p) > 0$.
  By Lemma~\ref{lem:cauchy} (iii), there exists $x \in S$ with $d(x, p) = d(S, p) > 0$, so that $x \ll p$.
  By Lemma~\ref{lem:cauchy} (i) applied to $S'$, there are three cases.

  First, suppose that $x \in I^+(S')$, and choose $y \in S'$ with $y \ll x$.
  Then the reverse triangle inequality shows $d(S', p) \geq d(y, p) \geq d(y, x) + d(x, p) \geq d(S, p)$.

  Next, suppose that $x \in S' \cup I^-(S')$ and $p \in J^+(S')$.
  Let $\gamma$ be a maximizing causal geodesic from $x$ to $p$ (Fact~\ref{fact:gh}).
  If $x \in S'$, we put $y \coloneqq x$; otherwise $\gamma$ meets $S'$ at some point $y$ by Lemma~\ref{lem:cauchy} (ii).
  Then we have
  \begin{equation}
    \begin{aligned}
      d(S, p) & = L_g(\gamma) = L_g(\gamma|_{x \to y}) + L_g(\gamma|_{y \to p}) \\
              & \leq d(x, y) + d(y, p) \leq d(S, S') + d(S', p) \leq R + d(S', p),
    \end{aligned}
  \end{equation}
  where we used Lemma~\ref{lem:max} for $d(S, S') \leq R$.

  Finally, suppose that $x \in S' \cup I^-(S')$ and $p \notin J^+(S')$.
  Then $p \in I^-(S')$ by Lemma~\ref{lem:cauchy} (i), so that there exists $y' \in S'$ with $p \ll y'$.
  Hence $d(S, p) = d(x, p) \leq d(x, p) + d(p, y') \leq d(x, y') \leq d(S, S') \leq R \leq d(S', p) + R$.
\end{proof}

\begin{lemma}\label{lem:ordered}
  Let $A, B \in \CM$ with $B \subset J^+(A)$. Then $d(B, A) = 0$, and hence $\dJ(A, B) = d(A, B)$.
\end{lemma}

\begin{proof}
  Suppose that $y \in B$, $x \in A$ and $d(y, x) > 0$, i.e., $y \ll x$.
  Since $y \in J^+(A)$, there exists $a \in A$ with $a \leq y$, so that $a \ll x$ with $a, x \in A$.
  This contradicts the achronality of $A$.
  The second claim follows from Lemma~\ref{lem:max}.
\end{proof}

We now show that $(\CM, \dJ)$ is proper,
using timelike geodesic completeness through Fact~\ref{fact:levels} instead of timelike Cauchy completeness (cf.~Fact~\ref{fact:LP} (ii)).
The graph argument below is the one used for local compactness in~\cite[Section 5.6]{LP};
the difference is that the parallel sets of a fixed Cauchy hypersurface give a compact region
containing a closed ball of any radius (Step~2), so that the argument yields the compactness of all closed balls directly.

\begin{proposition}\label{prop:proper}
  The metric space $(\CM, \dJ)$ is proper.
\end{proposition}

\begin{proof}
  \emph{Step 1: reduction to closed balls centered at $\Sigma$.}
  Fix an acausal compact Cauchy hypersurface $\Sigma$.
  Since every closed ball in $\CM$ is closed and contained in a closed ball centered at $\Sigma$,
  it suffices to show that every sequence in $\overline{B}(\Sigma, R) \coloneqq \{S \in \CM ; \dJ(S, \Sigma) \leq R\}$
  has a subsequence converging in $\dJ$ to an element of $\overline{B}(\Sigma, R)$.

  \emph{Step 2: construction of a compact set containing all elements of $\overline{B}(\Sigma, R)$.}
  Let $S \in \overline{B}(\Sigma, R)$ and $q \in S$.
  If $q \in J^+(\Sigma)$, then $d_\Sigma(q) = d(\Sigma, q) \leq d(\Sigma, S) \leq \dJ(\Sigma, S) \leq R$,
  and similarly $d_\Sigma(q) \geq -R$ if $q \in J^-(\Sigma)$.
  Hence $S \subset K_R \coloneqq d_\Sigma^{-1}([-R, R])$.
  The set $K_R$ is the union of $J^+(\Sigma) \cap J^-(\Sigma_R)$ and $J^-(\Sigma) \cap J^+(\Sigma_{-R})$,
  and it is compact by Fact~\ref{fact:levels} and its time dual;
  this is where timelike geodesic completeness is used.

  \emph{Step 3: graph representation of Cauchy hypersurfaces.}
  By~\cite{BS}, \cite[Theorem 4.120]{MiLR}, there exist a smooth Cauchy temporal function $\mathcal{T}$
  and a diffeomorphism $M \cong \R \times S'$ with $S'$ compact,
  in which $\mathcal{T}$ is the first coordinate and $\partial_{\mathcal{T}}$ is timelike.
  Since each integral curve $\R \times \{x\}$ is an inextendible timelike curve,
  each $S \in \CM$ is the graph $\{\mathcal{T} = f_S(x)\}$ of a function $f_S \colon S' \to \R$.

  \emph{Step 4: uniform Lipschitz estimate for graphs.}
  Choose $a < b$ with $K_R \subset \mathcal{T}^{-1}([a, b])$ and a Riemannian metric $h_0$ on $S'$.
  Since $\partial_{\mathcal{T}}$ is future-directed timelike and the set of $h_0$-unit vectors in $TS'$
  at points of $[a, b] \times S'$ is compact, there exists $\kappa > 0$ such that every vector
  $\lambda \partial_{\mathcal{T}} + w$ with $w \in TS'$ and $\lambda > \kappa \abs{w}_{h_0}$
  at a point of $[a, b] \times S'$ is future-directed timelike.
  We claim that if $S \subset K_R$, then $\abs{f_S(x) - f_S(y)} \leq \kappa \operatorname{dist}_{h_0}(x, y)$ for all $x, y \in S'$.
  Suppose otherwise. Then $x \neq y$, and we may assume $f_S(x) < f_S(y)$.
  Put $\ell \coloneqq \operatorname{dist}_{h_0}(x, y)$, and let $\sigma \colon [0, \ell] \to S'$ be a unit speed minimizing $h_0$-geodesic
  from $x$ to $y$. The curve
  \begin{equation}
    c(s) \coloneqq \Bigl( f_S(x) + \frac{s}{\ell} \bigl( f_S(y) - f_S(x) \bigr), \sigma(s) \Bigr), \qquad s \in [0, \ell],
  \end{equation}
  lies in $[a, b] \times S'$, since $f_S(x), f_S(y) \in [a, b]$ by $S \subset K_R \subset \mathcal{T}^{-1}([a, b])$.
  Its velocity is $\lambda \partial_{\mathcal{T}} + \dot\sigma(s)$ with $\lambda = (f_S(y) - f_S(x))/\ell > \kappa = \kappa \abs{\dot\sigma(s)}_{h_0}$,
  which is future-directed timelike.
  Thus $c$ is a timelike curve from $(f_S(x), x) \in S$ to $(f_S(y), y) \in S$, which contradicts the achronality of $S$.

  \emph{Step 5: taking a limit graph.}
  Let $(S_j)_j$ be a sequence in $\overline{B}(\Sigma, R)$ with graph functions $f_j$.
  The functions $f_j$ take values in $[a, b]$ and are $\kappa$-Lipschitz.
  By the Arzel\`a--Ascoli theorem, a subsequence, still denoted by $(f_j)_j$,
  converges uniformly to a $\kappa$-Lipschitz function $f$. Let $S$ be the graph of $f$.
  We show that $S \in \CM$ and $\dJ(S_j, S) \to 0$.

  \emph{Step 6: the limit graph is a Cauchy hypersurface.}
  The set $S$ is achronal. Indeed, if $p \ll q$ for $p, q \in S$,
  then the corresponding points of $S_j$ converge to $p$ and $q$,
  and since $\ll$ is open, they are chronologically related for large $j$.
  This contradicts the achronality of $S_j$.
  Moreover, $\mathcal{T}$ takes all real values along every inextendible timelike curve $\gamma$,
  so that $\mathcal{T}(\gamma) - f(\pr(\gamma))$ changes its sign, where $\pr \colon M \to S'$ is the projection.
  Thus $\gamma$ meets $S$, and it meets $S$ at most once by achronality.

  \emph{Step 7: $S_j$ converges to $S$ in $\dJ$.}
  Suppose not. By Lemma~\ref{lem:max}, there exist $\varepsilon > 0$ and, along a subsequence,
  $x_j \in S_j$ and $y_j \in S$ with $d(x_j, y_j) \geq \varepsilon$ or $d(y_j, x_j) \geq \varepsilon$.
  By the compactness of $[a, b] \times S'$, we may assume $x_j \to x$ and $y_j \to y$.
  We have $y \in S$, since $S$ is closed.
  Writing $x_j = (f_j(u_j), u_j)$ and $x = (t, u)$, we have $u_j \to u$ and
  $\abs{f_j(u_j) - f(u)} \leq \sup_{S'} \abs{f_j - f} + \abs{f(u_j) - f(u)} \to 0$
  by the uniform convergence, so that $t = f(u)$, i.e., $x \in S$.
  The continuity of $d$ gives $d(x, y) \geq \varepsilon$ or $d(y, x) \geq \varepsilon$, i.e., $x \ll y$ or $y \ll x$,
  which contradicts the achronality of $S$.

  Finally, since $\dJ(S, \Sigma) \leq \dJ(S, S_j) + \dJ(S_j, \Sigma) \leq \dJ(S, S_j) + R$ for all $j$,
  we obtain $S \in \overline{B}(\Sigma, R)$.
\end{proof}

Fix $\mathcal{T}$, $M \cong \R \times S'$ and the graph functions $f_S$ as in Step~3 of the proof of Proposition~\ref{prop:proper},
and put $H_t \coloneqq \mathcal{T}^{-1}(t)$ for $t \in \R$.
Fix a Riemannian metric $h_0$ on $S'$, let $\rho_0$ be the distance of the product metric $dt^2 + h_0$ on $\R \times S' \cong M$,
where $t$ is the coordinate of $\R$, which corresponds to $\mathcal{T}$,
and let $d_H$ be the Hausdorff distance between subsets of $(M, \rho_0)$.

\begin{lemma}\label{lem:topology}
  For every $S \in \CM$, the function $f_S$ is continuous.
  Moreover, for $S \in \CM$ and a sequence $(S_j)_j$ in $\CM$, the following are equivalent:
  \begin{enumerate}
    \item $\dJ(S_j, S) \to 0$;
    \item $\sup_{S'} \abs{f_{S_j} - f_S} \to 0$;
    \item $d_H(S_j, S) \to 0$.
  \end{enumerate}
  In particular, the map $S \mapsto f_S$ is a homeomorphism from $(\CM, \dJ)$
  onto its image in the space $C^0(S')$ of continuous functions with the uniform topology,
  and the topology of $\dJ$ coincides with the topology of $d_H$ on $\CM$.
\end{lemma}

\begin{proof}
  Each $S \in \CM$ is compact by~\cite[Corollary 14.32]{On},
  and the projection $S \to S'$ is a continuous bijection, since each integral curve $\R \times \{x\}$ meets $S$ exactly once.
  Hence it is a homeomorphism from a compact space onto a Hausdorff space, so that $f_S$ is continuous.
  We write $f_j \coloneqq f_{S_j}$ and $f \coloneqq f_S$.

  (ii) $\Rightarrow$ (i).
  Since $f$ is bounded and $f_j \to f$ uniformly, there exist $a < b$ such that $S$ and all $S_j$ lie in $[a, b] \times S'$.
  Then Step~7 of the proof of Proposition~\ref{prop:proper} applies without change,
  since it uses only the uniform convergence, the continuity of $f$, the closedness and achronality of $S$, and the compactness of $[a, b] \times S'$.

  (ii) $\Rightarrow$ (iii).
  For $x \in S'$, the integral curve of $\partial_{\mathcal{T}}$ between $(f_j(x), x)$ and $(f(x), x)$ has $\rho_0$-length $\abs{f_j(x) - f(x)}$.
  Hence every point of $S_j$ lies within $\rho_0$-distance $\sup_{S'} \abs{f_j - f}$ of $S$, and vice versa,
  so that $d_H(S_j, S) \leq \sup_{S'} \abs{f_j - f}$.

  (i) $\Rightarrow$ (ii).
  Suppose that (ii) fails.
  Then there exist $\varepsilon > 0$ and, along a subsequence, $u_j \in S'$ with $\abs{f_j(u_j) - f(u_j)} \geq \varepsilon$.
  Since $S'$ is compact, we may assume $u_j \to u$.
  Passing to a further subsequence and using the time duality, we may assume $f_j(u_j) \geq f(u_j) + \varepsilon$ for all $j$.
  Put $x_j \coloneqq (f_j(u_j), u_j) \in S_j$, $p_j \coloneqq (f(u_j), u_j) \in S$ and $q_j \coloneqq (f(u_j) + \varepsilon, u_j)$.
  Since $\partial_{\mathcal{T}}$ is future-directed timelike, we have $q_j \leq x_j$,
  and the reverse triangle inequality and Lemma~\ref{lem:max} give
  $\dJ(S, S_j) \geq d(p_j, x_j) \geq d(p_j, q_j)$.
  As $j \to \infty$, we have $p_j \to p \coloneqq (f(u), u)$ and $q_j \to q \coloneqq (f(u) + \varepsilon, u)$,
  and $d(p, q) > 0$ since $p \ll q$.
  The continuity of $d$ gives $\liminf_j \dJ(S, S_j) \geq d(p, q) > 0$, so that (i) fails.
  The idea of comparing points on the same integral curve of a timelike vector field
  is that of the estimate~(6) in the proof of~\cite[Proposition 5.23]{LP}, where a unit timelike vector field is used.

  (iii) $\Rightarrow$ (ii).
  Suppose that (ii) fails, and choose $\varepsilon$, $u_j$, $u$ and $x_j$ as in the proof of (i) $\Rightarrow$ (ii),
  where the time duality is not needed.
  If $d_H(S_j, S) \to 0$, then there exist $y_j \in S$ with $\rho_0(x_j, y_j) \to 0$.
  Since $S$ is compact, we may assume $y_j \to y \in S$, so that $x_j \to y$.
  Since $u_j \to u$, we have $y = (f(u), u)$, and hence $f_j(u_j) \to f(u)$.
  This contradicts $\abs{f_j(u_j) - f(u_j)} \geq \varepsilon$ and the continuity of $f$.
  Thus (iii) fails.

  The last statements follow, since all the topologies involved are metrizable,
  so that continuity can be checked by sequences.
\end{proof}

Lemma~\ref{lem:topology} uses neither timelike geodesic completeness nor any curvature condition,
and it shows that the topology of $\dJ$ depends only on the causal structure of $M$: 
if $\tilde g$ is a Lorentzian metric conformal to $g$, then $\tilde g$ and $g$ have the same Cauchy hypersurfaces and the same Cauchy temporal functions,
so that the metrics $\dJ$ of $g$ and $\tilde g$ define the same topology on $\CM$, although they are different in general.
Since the integral curves of $\partial_{\mathcal{T}}$ are timelike, Lemma~\ref{lem:cauchy} (i) gives $I^-(S) = \{\mathcal{T} < f_S\}$ for every $S \in \CM$.

We extend the signed distance of Section~\ref{subsec:busemann} to all $S \in \CM$ by
\begin{equation}
  d_S(p) \coloneqq d(S, p) - d(p, S), \qquad p \in M.
\end{equation}
If $p \in J^+(S)$, then $d(p, S) = 0$, since $s' \leq p \ll s$ with $s, s' \in S$ would contradict the achronality of $S$; hence $d_S = d(S, \cdot)$ on $J^+(S)$, and similarly $d_S = -d(\cdot, S)$ on $J^-(S)$.
The level sets $S_c$ defined below extend the notation $\Sigma_c$ of Section~\ref{subsec:busemann}.
The function $d_S$ and its level sets $S_c$ are the function $\tau_S$ and the parallel sets of $S$
considered in~\cite[Section 5.5]{LP};
for timelike Cauchy complete spacetimes, the fact that the level sets $S_c$ are Cauchy hypersurfaces is contained in~\cite[Lemma 5.16]{LP}.
\begin{definition}\label{def:monotone}
  For $T, T' \in \CM$, we write $T \preceq T'$ if $T \subset J^-(T')$.
  A family $(S_t)_{t \in \R}$ in $\CM$ is called \emph{causally monotone}
  if $S_t \subset J^+(S_s)$ for all $s < t$, or equivalently (Lemma~\ref{lem:lattice} (i) below), if $S_s \preceq S_t$ for all $s < t$.
\end{definition}

\begin{lemma}\label{lem:signed}
  Let $S \in \CM$. Then the following hold:
  \begin{enumerate}
    \item $d_S$ is continuous, $d_S(p) < d_S(q)$ for all $p \ll q$,
          and $\abs{d_S - d_\Sigma} \leq 2\dJ(S, \Sigma)$ for every acausal compact Cauchy hypersurface $\Sigma$;
    \item for every $c \in \R$, the level set $S_c \coloneqq d_S^{-1}(c)$ is a Cauchy hypersurface, $S_0 = S$, and
          \begin{equation}
            \{d_S \leq c\} = J^-(S_c), \qquad \{d_S \geq c\} = J^+(S_c).
          \end{equation}
  \end{enumerate}
\end{lemma}

\begin{proof}
  (i) The continuity follows from Lemma~\ref{lem:cauchy} (iii) and its time dual,
  and the last inequality follows from Lemma~\ref{lem:lipschitz}.
  Let $p \ll q$. By Lemma~\ref{lem:cauchy} (i), there are three cases.
  If $p \in J^+(S)$, then $q \in I^+(S)$.
  If $d(S, p) > 0$, choose $x \in S$ with $d(x, p) = d(S, p)$ by Lemma~\ref{lem:cauchy} (iii);
  then $d(S, q) \geq d(x, q) \geq d(x, p) + d(p, q) > d(S, p)$.
  If $d(S, p) = 0$, choose $s \in S$ with $s \leq p$; then $s \ll q$ and $d(S, q) \geq d(s, q) > 0$.
  In both cases, $d_S(q) > d_S(p)$.
  If $q \in J^-(S)$, the time-dual argument applies.
  Otherwise, $p \in I^-(S)$ and $q \in I^+(S)$, so that $d_S(p) < 0 < d_S(q)$.

  (ii) Fix an acausal compact Cauchy hypersurface $\Sigma$ (Section~\ref{sec:prelim}).
  Let $\gamma$ be a future-directed inextendible timelike curve.
  Since $d_\Sigma$ is a Cauchy time function (Fact~\ref{fact:levels}), $d_\Sigma \circ \gamma$ takes all real values,
  and so does the continuous function $d_S \circ \gamma$ by (i).
  Since $d_S \circ \gamma$ is strictly increasing by (i), $\gamma$ meets $S_c$ exactly once.
  Thus $S_c \in \CM$.
  By the remark before Definition~\ref{def:monotone} and the achronality of $S$,
  we have $d_S > 0$ on $I^+(S)$, $d_S < 0$ on $I^-(S)$ and $d_S = 0$ on $S$, so that $S_0 = S$ by Lemma~\ref{lem:cauchy} (i).
  If $d_S(p) < c$, then a future-directed inextendible timelike curve from $p$ meets $S_c$ at a point $q$ with $p \ll q$,
  so that $p \in I^-(S_c)$; conversely, $d_S < c$ on $I^-(S_c)$ by (i).
  Hence $\{d_S < c\} = I^-(S_c)$, and Lemma~\ref{lem:cauchy} (i) gives $\{d_S \leq c\} = I^-(S_c) \cup S_c = J^-(S_c)$.
  The last equality is the time-dual statement.
\end{proof}

\begin{lemma}\label{lem:lattice}
  The following hold:
  \begin{enumerate}
    \item for $T, T' \in \CM$, the conditions $T \preceq T'$, $T' \subset J^+(T)$ and $I^-(T) \subset I^-(T')$ are equivalent.
          Moreover, $\overline{I^-(T)} = J^-(T)$ and $\partial I^-(T) = T$.
          In particular, $\preceq$ is a partial order on $\CM$;
    \item if $A \subset M$ is an open set with $I^-(A) \subset A$ such that every inextendible timelike curve
          meets both $A$ and $M \setminus A$, then $\partial A \in \CM$ and $I^-(\partial A) = A$;
    \item every finite non-empty family $T_1, \dots, T_k \in \CM$ has a least upper bound and a greatest lower bound
          with respect to $\preceq$, given by
          \begin{equation}
            \bigvee_{i} T_i = \partial \Bigl( \bigcup_i I^-(T_i) \Bigr), \qquad
            \bigwedge_{i} T_i = \partial \Bigl( \bigcap_i I^-(T_i) \Bigr).
          \end{equation}
  \end{enumerate}
\end{lemma}

\begin{proof}
  (i) By Lemma~\ref{lem:cauchy} (i), $M = I^-(T) \sqcup T \sqcup I^+(T)$ with $I^\pm(T)$ open,
  and every point of $T$ is a limit of points of $I^-(T)$ and of $I^+(T)$. Hence
  \begin{equation}
    J^-(T) = M \setminus I^+(T) = \overline{I^-(T)}, \qquad I^-(T) = \operatorname{Int} J^-(T), \qquad \partial I^-(T) = T,
  \end{equation}
  and similarly $J^+(T) = M \setminus I^-(T)$. Therefore
  \begin{equation}
    \begin{aligned}
      T \subset J^-(T') & \iff J^-(T) \subset J^-(T') \iff I^-(T) \subset I^-(T') \\
                        & \iff J^+(T') \subset J^+(T) \iff T' \subset J^+(T),
    \end{aligned}
  \end{equation}
  where the first and last equivalences follow from the transitivity of $\leq$,
  the second one from taking closures and interiors, and the third one from taking complements.
  The relation $\preceq$ is clearly reflexive and transitive, and it is antisymmetric
  since $T \preceq T' \preceq T$ gives $I^-(T) = I^-(T')$, and hence $T = \partial I^-(T) = \partial I^-(T') = T'$.

  (ii) Since $A$ is open, $A = I^-(A)$ is a past set,
  so $\partial A$ is achronal~\cite[Remark 2.85 and Theorem 2.87 (iii)]{MiLR}.
  Let $\gamma$ be an inextendible future-directed timelike curve.
  Since $A$ is open and $I^-(A) \subset A$, the set of parameters $s$ with $\gamma(s) \in A$ is an open initial segment,
  which is neither empty nor the whole parameter interval by assumption.
  Hence it is of the form $\{s < s_0\}$, and $\gamma$ meets $\partial A$ at $\gamma(s_0)$, exactly once by achronality.
  Thus $\partial A \in \CM$ and $A \subset I^-(\partial A)$;
  conversely, $I^-(\partial A) \subset I^-(\overline{A}) = I^-(A) \subset A$, since $I^+(p)$ is open for every $p \in M$.

  (iii) The set $A \coloneqq \bigcup_i I^-(T_i)$ is open and satisfies $I^-(A) \subset A$.
  An inextendible timelike curve meets $A$ before it meets $T_1$,
  and meets $\bigcap_i I^+(T_i) \subset M \setminus A$ after it has met all $T_i$.
  By (ii), $X \coloneqq \partial A \in \CM$ and $I^-(X) = A$,
  and by (i), $X \preceq T$ if and only if $A \subset I^-(T)$, i.e., $T_i \preceq T$ for all $i$.
  Thus the upper bounds of $T_1, \dots, T_k$ are exactly the $T \in \CM$ with $X \preceq T$,
  so $X$ is the least upper bound.
  The meet is obtained in the same way from $\bigcap_i I^-(T_i)$.
\end{proof}

\section{Injectivity of $(\CM, \dJ)$}\label{sec:injective}

In this section, we show that closed balls in $(\CM, \dJ)$ are order intervals,
and deduce that $(\CM, \dJ)$ is an injective metric space.

For $L, U \in \CM$, we write $[L, U] \coloneqq \{T \in \CM ; L \preceq T \preceq U\}$.

\begin{proposition}\label{prop:ball}
  For every $S \in \CM$ and $r \geq 0$, we have $\overline{B}(S, r) = [S_{-r}, S_r]$.
\end{proposition}

\begin{proof}
  Let $y \in M$. By Lemma~\ref{lem:max} and the remark after the definition of $d_S$,
  we have $\sup_{x \in S} (d(x, y) + d(y, x)) = \max\{d(S, y), d(y, S)\} = \abs{d_S(y)}$,
  and hence $\dJ(S, T) = \sup_{y \in T} \abs{d_S(y)}$ for $T \in \CM$.
  By Lemma~\ref{lem:signed} (ii) and Lemma~\ref{lem:lattice} (i), we obtain
  \begin{equation}
    \dJ(S, T) \leq r \iff T \subset \{-r \leq d_S \leq r\} = J^+(S_{-r}) \cap J^-(S_r) \iff S_{-r} \preceq T \preceq S_r.
  \end{equation}
\end{proof}

Following~\cite{La}, a metric space $(X, \rho)$ is \emph{injective} if, for every metric space $Y$, every $A \subset Y$,
and every $1$-Lipschitz map $f \colon A \to X$, there exists a $1$-Lipschitz extension $Y \to X$ of $f$,
and \emph{hyperconvex} if, for every family $(x_i, r_i)_{i \in I}$ in $X \times \R$
with $\rho(x_i, x_j) \leq r_i + r_j$ for all $i, j \in I$, we have $\bigcap_{i \in I} \overline{B}(x_i, r_i) \neq \emptyset$
(note that $r_i \geq 0$ by taking $i = j$, and that the intersection over the empty family is $X$).
These two conditions are equivalent~\cite{AP}; see also~\cite[Proposition 2.3]{La}.

\begin{theorem}\label{thm:injective}
  Let $(S_i, r_i)_{i \in I}$ be a family in $\CM \times \R$ with $\dJ(S_i, S_j) \leq r_i + r_j$ for all $i, j \in I$.
  Then $\bigcap_{i \in I} \overline{B}(S_i, r_i) \neq \emptyset$.
  If $I$ is finite and non-empty, then
  \begin{equation}\label{eq:helly}
    \bigcap_{i \in I} \overline{B}(S_i, r_i) = \Bigl[ \bigvee_{i \in I} (S_i)_{-r_i}, \ \bigwedge_{i \in I} (S_i)_{r_i} \Bigr] \neq \emptyset.
  \end{equation}
  Consequently, $(\CM, \dJ)$ is a hyperconvex, and hence injective, metric space.
  In particular, it is a proper geodesic metric space.
\end{theorem}

\begin{proof}
  \emph{Step 1: we show that $(S_i)_{-r_i} \preceq (S_j)_{r_j}$ for all $i, j \in I$.}
  This step is a variant of the argument in the proof of~\cite[Lemma 5.19]{LP}.
  Suppose not. By Lemma~\ref{lem:signed} (ii) and Lemma~\ref{lem:lattice} (i),
  there exists $p \in (S_i)_{-r_i}$ with $d_{S_j}(p) > r_j$.
  Then $d(S_j, p) > r_j$, and $d_{S_i}(p) = -r_i$ gives $d(p, S_i) = r_i$.
  By Lemma~\ref{lem:cauchy} (iii) and its time dual, there exist $x \in S_j$ and $y \in S_i$
  with $d(x, p) = d(S_j, p) > r_j$ and $d(p, y) = d(p, S_i) = r_i$,
  so that $x \ll p$, and $p \ll y$ if $r_i > 0$; if $r_i = 0$, then $p \in S_i$ and we take $y = p$ instead.
  Then $x \ll p \leq y$, and the reverse triangle inequality gives
  \begin{equation}
    \dJ(S_i, S_j) \geq d(x, y) \geq d(x, p) + d(p, y) > r_j + r_i,
  \end{equation}
  which contradicts the assumption $\dJ(S_i, S_j) \leq r_i + r_j$.

  \emph{Step 2: the finite family case.}
  Suppose that $I$ is finite and non-empty, and let $L \coloneqq \bigvee_i (S_i)_{-r_i}$ and $U \coloneqq \bigwedge_i (S_i)_{r_i}$,
  which exist by Lemma~\ref{lem:lattice} (iii).
  By Proposition~\ref{prop:ball}, $T \in \bigcap_i \overline{B}(S_i, r_i)$ if and only if
  $(S_i)_{-r_i} \preceq T \preceq (S_i)_{r_i}$ for all $i$, i.e., $L \preceq T \preceq U$.
  Hence $\bigcap_i \overline{B}(S_i, r_i) = [L, U]$, and it remains to show that $L \preceq U$, i.e., that $[L, U] \neq \emptyset$.
  By the first step, each $(S_j)_{r_j}$ is an upper bound of the family $((S_i)_{-r_i})_i$,
  so that $L \preceq (S_j)_{r_j}$ for all $j$, and hence $L \preceq U$ since $U$ is the greatest lower bound of $((S_j)_{r_j})_j$.
  This proves~\eqref{eq:helly}.

  \emph{Step 3: the arbitrary family case.}
  We may assume $I \neq \emptyset$, since the intersection over the empty family is $\CM \neq \emptyset$.
  Fix $i_0 \in I$ and let $B \coloneqq \overline{B}(S_{i_0}, r_{i_0})$, which is compact by Proposition~\ref{prop:proper}.
  The family of sets $B \cap \overline{B}(S_i, r_i)$, $i \in I$, consists of closed subsets of $B$,
  and it has the finite intersection property by~\eqref{eq:helly}.
  Hence its intersection, which is $\bigcap_{i \in I} \overline{B}(S_i, r_i)$, is non-empty.

  \emph{Step 4: consequences.}
  Thus $(\CM, \dJ)$ is hyperconvex, and hence injective.
  It is proper by Proposition~\ref{prop:proper},
  and geodesic since every injective metric space is~\cite[Section 2]{La}.
\end{proof}

\section{Distance curves}\label{sec:parallel}

In this section, we study the curve in $\CM$ formed by the parallel sets of a fixed acausal compact
Cauchy hypersurface $\Sigma$, and reformulate the criterion of~\cite{Mi} (Fact~\ref{fact:defect})
in terms of the metric $\dJ$.
Throughout this section, $\Sigma$ denotes an acausal compact Cauchy hypersurface.

\subsection{The distance curve and its defect}

\begin{lemma}\label{lem:monotone-family}
  If $c < c'$, then $\Sigma_{c'} \subset I^+(\Sigma_c)$.
  In particular, the family $(\Sigma_t)_{t \in \R}$ is causally monotone,
  and $\dJ(\Sigma_c, \Sigma_{c'}) = d(\Sigma_c, \Sigma_{c'})$.
\end{lemma}

\begin{proof}
  Let $y \in \Sigma_{c'}$. An inextendible timelike curve through $y$ meets the Cauchy hypersurface $\Sigma_c$
  at a point $x$ (Fact~\ref{fact:levels}).
  Since $d_\Sigma$ is a time function and $d_\Sigma(x) = c < c' = d_\Sigma(y)$,
  the point $x$ must precede $y$ on the curve, so that $x \ll y$.
  The last claim follows from Lemma~\ref{lem:ordered}.
\end{proof}

For $a, b \geq 0$, we define the \emph{defect}
\begin{equation}
  D(a, b) \coloneqq a + b - d(\Sigma_{-a}, \Sigma_b).
\end{equation}
We also set $\Delta_k \coloneqq D(k, k) = 2k - d(\Sigma_{-k}, \Sigma_k)$ for $k \in \N$; this is the sequence considered in~\cite{Mi} (see~Fact~\ref{fact:defect}).

\begin{lemma}\label{lem:defect}
  For all $a, b \geq 0$, we have $0 \leq D(a, b) \leq a + b$,
  and $D$ is non-decreasing in each variable.
\end{lemma}

\begin{proof}
  The lower bound is $d(\Sigma_{-a}, \Sigma_b) \leq a + b$, which is Fact~\ref{fact:lengths},
  and the upper bound follows from $d \geq 0$.

  We first show the monotonicity in $b$. Let $0 \leq b \leq b'$. We show
  \begin{equation}\label{eq:step}
    d(\Sigma_{-a}, \Sigma_{b'}) \leq d(\Sigma_{-a}, \Sigma_b) + (b' - b),
  \end{equation}
  which is equivalent to $D(a, b) \leq D(a, b')$.
  If $b > 0$, this follows from Fact~\ref{fact:lengths} by taking the supremum over $y \in \Sigma_{-a}$.
  If $b = 0$, let $x \in \Sigma_{-a}$ and $z \in \Sigma_{b'}$ with $x \leq z$.
  Every causal curve $\gamma$ from $x$ to $z$ meets $\Sigma$ at some point $y$, so that
  $L_g(\gamma) \leq d(x, y) + d(y, z) \leq d(\Sigma_{-a}, \Sigma) + d(\Sigma, z) = d(\Sigma_{-a}, \Sigma) + b'$,
  and taking suprema gives~\eqref{eq:step}.
  The monotonicity in $a$ follows in the same way from the time-dual statements.
\end{proof}

\begin{definition}\label{def:parallel}
  The \emph{distance curve} of $\Sigma$ is the curve
  \begin{equation}
    c_\Sigma \colon \R \to \CM, \qquad c_\Sigma(t) \coloneqq \Sigma_t.
  \end{equation}
  We define
  \begin{equation}
    D_\infty(\Sigma) \coloneqq \sup_{a, b \geq 0} D(a, b) \in [0, \infty].
  \end{equation}
\end{definition}

The sets $\Sigma_t$ are Lorentzian analogues of the equidistant hypersurfaces of a hypersurface in Riemannian geometry,
and they are in general only Lipschitz.

\begin{proposition}\label{prop:independent}
  Let $\Sigma'$ be another acausal compact Cauchy hypersurface.
  Then the following hold:
  \begin{enumerate}
    \item $\dJ(\Sigma_t, \Sigma'_t) \leq \dJ(\Sigma, \Sigma')$ for every $t \in \R$;
    \item $\dom b^\pm_\Sigma = \dom b^\pm_{\Sigma'}$, and $\abs{b^\pm_\Sigma - b^\pm_{\Sigma'}} \leq \dJ(\Sigma, \Sigma')$ on this common domain.
  \end{enumerate}
\end{proposition}

\begin{proof}
  (ii) is~\cite[Theorem 3.9]{Mi}, since the constant $R(\Sigma, \Sigma')$ there coincides with $\dJ(\Sigma, \Sigma')$ by Lemma~\ref{lem:max}.
  (i) is shown in the proof of~\cite[Theorem 3.9]{Mi} for $\abs{t} > \dJ(\Sigma, \Sigma')$; we prove it for all $t$.

  (i) Put $R \coloneqq \dJ(\Sigma, \Sigma')$. The case $t = 0$ is trivial.
  Let $t > 0$, $p \in \Sigma_t$ and $q \in \Sigma'_t$. We show $d(p, q) \leq R$.
  Then $d(q, p) \leq R$ by exchanging the roles of $\Sigma$ and $\Sigma'$, and $\dJ(\Sigma_t, \Sigma'_t) \leq R$ follows from Lemma~\ref{lem:max}.
  If $p \not\leq q$, then $d(p, q) = 0$.
  If $p \leq q$, choose $x \in \Sigma$ with $d(x, p) = d(\Sigma, p) = t > 0$ (Lemma~\ref{lem:cauchy} (iii)).
  Then $x \ll p \leq q$, and the reverse triangle inequality gives $d(\Sigma, q) \geq d(x, q) \geq t + d(p, q)$.
  On the other hand, Lemma~\ref{lem:lipschitz} gives $d(\Sigma, q) \leq d(\Sigma', q) + R = t + R$.
  Hence $d(p, q) \leq R$.
  The case $t < 0$ is the time-dual statement.
\end{proof}

The following proposition shows that $c_\Sigma$ moves at unit speed on each side of $\Sigma$,
and that $D_\infty(\Sigma)$ measures how far $c_\Sigma$ is from being a geodesic line.
Across $\Sigma$, a causal curve from $\Sigma_{-a}$ to $\Sigma_b$ crosses $\Sigma$ at a single point,
and since $\Sigma$ is a hypersurface rather than a point, the points of $\Sigma$ realizing the distance from $\Sigma_{-a}$,
i.e., the points $y \in \Sigma$ with $d(\Sigma_{-a}, y) = a$,
and those realizing the distance to $\Sigma_b$ are different in general; this is why $D(a, b) > 0$ in general.

\begin{proposition}\label{prop:parallel}
  The following hold:
  \begin{enumerate}
    \item if $s, t \geq 0$ or $s, t \leq 0$, then $\dJ(\Sigma_s, \Sigma_t) = \abs{t - s}$.
          In particular, $c_\Sigma|_{[0, \infty)}$ and $c_\Sigma|_{(-\infty, 0]}$ are geodesic rays;
    \item for $a, b \geq 0$, we have $\dJ(\Sigma_{-a}, \Sigma_b) = a + b - D(a, b)$, and
          \begin{equation}
            D_\infty(\Sigma) = \lim_{k \to \infty} \Delta_k \in [0, \infty];
          \end{equation}
    \item for all $s, t \in \R$, we have
          \begin{equation}
            \abs{t - s} - D_\infty(\Sigma) \leq \dJ(\Sigma_s, \Sigma_t) \leq \abs{t - s},
          \end{equation}
          and $D_\infty(\Sigma)$ is the smallest constant with this property.
          In particular, $c_\Sigma$ is a $(1, \varepsilon)$-quasi-geodesic line for some $\varepsilon \geq 0$
          if and only if $D_\infty(\Sigma) < \infty$,
          and $c_\Sigma$ is a geodesic line if and only if $D_\infty(\Sigma) = 0$;
    \item if $\Sigma'$ is another acausal compact Cauchy hypersurface, then, in $[0, \infty]$,
          \begin{equation}
            D_\infty(\Sigma') \leq D_\infty(\Sigma) + 2\dJ(\Sigma, \Sigma'), \qquad
            D_\infty(\Sigma) \leq D_\infty(\Sigma') + 2\dJ(\Sigma, \Sigma').
          \end{equation}
          In particular, $D_\infty(\Sigma) < \infty$ if and only if $D_\infty(\Sigma') < \infty$.
  \end{enumerate}
\end{proposition}

\begin{proof}
  (i) For timelike Cauchy complete spacetimes, this is~\cite[Lemma 5.22]{LP}; we give a short proof.
  By the symmetry of $\dJ$, we may assume $s < t$. Let $0 \leq s < t$.
  By Lemma~\ref{lem:monotone-family}, we have $\dJ(\Sigma_s, \Sigma_t) = d(\Sigma_s, \Sigma_t)$.
  By Fact~\ref{fact:levels}, or by the definition of $\Sigma_t$ if $s = 0$, we have $\Sigma_t = (d^+_{\Sigma_s})^{-1}(t - s)$,
  so that every $y \in \Sigma_t$ satisfies $d(\Sigma_s, y) = t - s$, and hence $d(\Sigma_s, \Sigma_t) = t - s$.
  The case $s, t \leq 0$ is the time-dual statement.

  (ii) If $a = b = 0$, both sides of the first formula vanish by achronality.
  Otherwise, it follows from Lemma~\ref{lem:monotone-family} and the definition of $D$.
  For the second formula, since $D$ is non-decreasing in each variable by Lemma~\ref{lem:defect},
  we have $D(a, b) \leq D(k, k) = \Delta_k$ for $k \in \N$ with $k \geq \max\{a, b\}$, and $(\Delta_k)_k$ is non-decreasing.
  Hence $\sup_{a, b \geq 0} D(a, b) = \sup_k \Delta_k = \lim_{k \to \infty} \Delta_k$ in $[0, \infty]$.

  (iii) If $s$ and $t$ have the same sign, this is (i).
  If $s = -a \leq 0 \leq t = b$, then $\abs{t - s} = a + b$, and (ii) gives
  $\dJ(\Sigma_s, \Sigma_t) = a + b - D(a, b) \in [a + b - D_\infty(\Sigma), a + b]$.
  The case $t \leq 0 \leq s$ follows from the symmetry of $\dJ$.
  The optimality of $D_\infty(\Sigma)$ follows from the definition of the supremum,
  and the remaining statements follow immediately.

  (iv) Let $R \coloneqq \dJ(\Sigma, \Sigma')$ and $k \in \N$, and let $\Delta'_k \coloneqq 2k - d(\Sigma'_{-k}, \Sigma'_k)$.
  By the triangle inequality for $\dJ$ and Proposition~\ref{prop:independent} (i), we have
  \begin{equation}
    \dJ(\Sigma_{-k}, \Sigma_k) \leq \dJ(\Sigma_{-k}, \Sigma'_{-k}) + \dJ(\Sigma'_{-k}, \Sigma'_k) + \dJ(\Sigma'_k, \Sigma_k)
    \leq \dJ(\Sigma'_{-k}, \Sigma'_k) + 2R.
  \end{equation}
  By Lemma~\ref{lem:monotone-family}, both sides are the corresponding values of $d$,
  so that $2k - \Delta_k \leq 2k - \Delta'_k + 2R$, i.e., $\Delta'_k \leq \Delta_k + 2R$.
  Letting $k \to \infty$ and using (ii), we obtain $D_\infty(\Sigma') \leq D_\infty(\Sigma) + 2R$.
  Exchanging $\Sigma$ and $\Sigma'$ gives the other inequality.
\end{proof}

Combined with the results of~\cite{Mi} on Cauchy Busemann functions, Proposition~\ref{prop:parallel} gives the following.

\begin{proposition}\label{prop:minguzzi}
  We have $D_\infty(\Sigma) = \min_\Sigma (b^+ - b^-) = \inf_M (b^+ - b^-)$.
  In particular, $c_\Sigma$ is a $(1, \varepsilon)$-quasi-geodesic line for some $\varepsilon \geq 0$
  if and only if $\dom b^+$ and $\dom b^-$ intersect.
\end{proposition}

\begin{proof}
  This follows from Proposition~\ref{prop:parallel} (ii), (iii) and Fact~\ref{fact:defect}.
\end{proof}

Under the timelike convergence condition, Fact~\ref{fact:defect} thus shows that $(M, g)$ splits if and only if
$c_\Sigma$ is a $(1, \varepsilon)$-quasi-geodesic line for some $\varepsilon \geq 0$.

The following lemma bounds the defect by the smaller of the two parameters.
Together with Lemma~\ref{lem:gromov}, which computes Gromov products along the distance curve,
it will be used in Section~\ref{sec:main} to relate Gromov hyperbolicity with the defect.

\begin{lemma}\label{lem:lower}
  For all $a, b \geq 0$, we have
  \begin{equation}
    \max\{a, b\} \leq \dJ(\Sigma_{-a}, \Sigma_b) \leq a + b.
  \end{equation}
  Equivalently, $D(a, b) \leq \min\{a, b\}$.
\end{lemma}

\begin{proof}
  By Proposition~\ref{prop:parallel} (ii), we have $\dJ(\Sigma_{-a}, \Sigma_b) = d(\Sigma_{-a}, \Sigma_b) = a + b - D(a, b)$,
  and the upper bound follows from Lemma~\ref{lem:defect}.
  We show $d(\Sigma_{-a}, \Sigma_b) \geq b$. We may assume $b > 0$.
  Let $y \in \Sigma_b$, and choose $x \in \Sigma$ with $d(x, y) = d(\Sigma, y) = b$ (Lemma~\ref{lem:cauchy} (iii)).
  If $a = 0$, then $d(\Sigma_{-a}, \Sigma_b) \geq d(x, y) = b$.
  If $a > 0$, then $\Sigma \subset I^+(\Sigma_{-a})$ by Lemma~\ref{lem:monotone-family},
  so that there exists $w \in \Sigma_{-a}$ with $w \ll x$.
  The reverse triangle inequality gives $d(\Sigma_{-a}, \Sigma_b) \geq d(w, y) \geq d(w, x) + d(x, y) \geq b$.
  The inequality $d(\Sigma_{-a}, \Sigma_b) \geq a$ is the time-dual statement.
\end{proof}

\begin{lemma}\label{lem:gromov}
  For $0 \leq r \leq s$, the Gromov products in $(\CM, \dJ)$ satisfy
  \begin{equation}
    (\Sigma_r \cdot \Sigma_s)_\Sigma = (\Sigma_{-r} \cdot \Sigma_{-s})_\Sigma = r, \qquad
    (\Sigma_r \cdot \Sigma_{-r})_\Sigma = \frac{1}{2} D(r, r) \leq \frac{r}{2}.
  \end{equation}
\end{lemma}

\begin{proof}
  By Proposition~\ref{prop:parallel} (i), we have $\dJ(\Sigma, \Sigma_{\pm r}) = r$, $\dJ(\Sigma, \Sigma_{\pm s}) = s$
  and $\dJ(\Sigma_{\pm r}, \Sigma_{\pm s}) = s - r$, which gives the first two equalities.
  By Proposition~\ref{prop:parallel} (ii), we have $\dJ(\Sigma_{-r}, \Sigma_r) = 2r - D(r, r)$,
  so that $(\Sigma_r \cdot \Sigma_{-r})_\Sigma = D(r, r)/2$, and $D(r, r) \leq r$ by Lemma~\ref{lem:lower}.
\end{proof}

\section{From quasi-geodesic lines to timelike lines}\label{sec:line}

In this section, we prove Theorem~\ref{thm:quasi}.
The quasi-geodesic line in Theorem~\ref{thm:quasi} is arbitrary: it need not be a distance curve,
and it need not even be causally monotone (Definition~\ref{def:monotone}).
Hence the arguments of Section~\ref{sec:parallel}, which use the parallel sets of a fixed Cauchy hypersurface, are not available.
Instead, we show that the distance from $S_0$ along a limit $\gamma$ of maximizing segments between $S_{-k}$ and $S_k$
is almost attained at $\gamma(0)$, which forces $\gamma$ to be timelike.
For $\varepsilon > 0$, the future half of $\gamma$ is a ray from $S_0$ only up to an additive error,
so that the fact that rays from a compact Cauchy hypersurface are timelike~\cite[Theorem 2.3]{Mi}
does not apply directly; Step~8 of the proof adapts its argument.

\begin{theorem}\label{thm:quasi}
  Let $(M, g)$ be a timelike geodesically complete, spatially compact, globally hyperbolic spacetime,
  and let $(S_t)_{t \in \R}$ be a $(1, \varepsilon)$-quasi-geodesic line in $(\CM, \dJ)$.
  Then there exists a complete timelike line $\gamma$ passing through a point $z \in S_0$.
  Moreover, if $\gamma$ is parametrized by proper time with $\gamma(0) = z$, then
  \begin{equation}
    d(S_0, \gamma(t)) \leq t + 3\varepsilon \quad (t \geq 0), \qquad
    d(\gamma(t), S_0) \leq \abs{t} + 3\varepsilon \quad (t \leq 0).
  \end{equation}
\end{theorem}

\begin{proof}
  We only use $S_{-k}$, $S_0$ and $S_k$ for integers $k > 2\varepsilon$.
  Applying Definition~\ref{def:quasi} with $(s, t) = (-k, k)$, $(0, k)$ and $(0, -k)$, we have
  \begin{equation}
    2k - \varepsilon \leq \dJ(S_{-k}, S_k) \leq 2k + \varepsilon, \qquad \dJ(S_0, S_k) \leq k + \varepsilon, \qquad \dJ(S_0, S_{-k}) \leq k + \varepsilon.
  \end{equation}

  \emph{Step 1: maximizing segments.}
  By Lemma~\ref{lem:max}, we have
  \begin{equation}
    L_k \coloneqq \dJ(S_{-k}, S_k) = \max\{d(S_{-k}, S_k), d(S_k, S_{-k})\} \in [2k - \varepsilon, 2k + \varepsilon].
  \end{equation}
  By the compactness of $S_{\pm k}$ and the continuity of $d$ (Fact~\ref{fact:gh}), the maximum is attained:
  there exist $x_k \in P_k$ and $y_k \in F_k$ with $L_k = d(x_k, y_k)$, where $\{P_k, F_k\} = \{S_{-k}, S_k\}$.
  It does not matter which of the two hypersurfaces contains the past endpoint $x_k$ (see Figure~\ref{fig:quasi} (a)).
  In either case, we have
  \begin{equation}
    \dJ(S_0, P_k) \leq k + \varepsilon, \qquad \dJ(S_0, F_k) \leq k + \varepsilon.
  \end{equation}
  The theorem of Avez and Seifert (Fact~\ref{fact:gh}) gives a maximizing causal geodesic $\gamma_k$
  from $x_k$ to $y_k$ of length $L_k$, which is timelike since $L_k > 0$.

  \emph{Step 2: $\gamma_k$ meets $S_0$ between its endpoints.}
  We extend $\gamma_k$ to a complete timelike geodesic.
  It meets $S_0$ at exactly one point $z_k$.
  Suppose that $z_k$ does not lie on the segment from $x_k$ to $y_k$.
  Then, since the extended curve is a timelike geodesic, either $z_k \ll x_k$ or $y_k \ll z_k$.
  If $z_k \ll x_k$, then the reverse triangle inequality gives
  \begin{equation}
    d(z_k, y_k) \geq d(z_k, x_k) + d(x_k, y_k) \geq L_k \geq 2k - \varepsilon > k + \varepsilon,
  \end{equation}
  since $k > 2\varepsilon$. Thus $d(S_0, F_k) \geq d(z_k, y_k) > k + \varepsilon$.
  However, $d(S_0, F_k) \leq \dJ(S_0, F_k) \leq k + \varepsilon$ by Lemma~\ref{lem:max} and Step~1, a contradiction.
  If $y_k \ll z_k$, then $d(x_k, z_k) \geq d(x_k, y_k) + d(y_k, z_k) \geq L_k > k + \varepsilon$,
  so that $d(P_k, S_0) > k + \varepsilon$.
  However, $d(P_k, S_0) \leq \dJ(P_k, S_0) = \dJ(S_0, P_k) \leq k + \varepsilon$, which is again a contradiction.
  Hence $z_k$ lies on $\gamma_k$.

  \begin{figure}[t]
    \centering
    \usetikzlibrary{intersections, decorations.pathreplacing, calc}
    \begin{minipage}[b]{0.5\textwidth}
      \centering
    \begin{tikzpicture}[scale=0.42, every node/.style={font=\footnotesize}]
      \draw[thick, name path=S0] plot[smooth, tension=0.7] coordinates {(0,-0.3) (2,0.3) (5,0) (8,0.15) (12,-0.05)};
      \draw[thick, blue, name path=Sk] plot[smooth, tension=0.7] coordinates {(0,-0.6) (1.2,-0.1) (3,1.1) (6,2.3) (9,2.1) (12,1.6)};
      \draw[thick, red, name path=Smk] plot[smooth, tension=0.7] coordinates {(0,0.5) (1.2,-0.1) (3,-1.2) (6,-2.3) (9,-2.1) (12,-1.6)};
      \node[right] at (12,-0.05) {$S_0$};
      \node[right, blue] at (12,1.6) {$S_k$};
      \node[right, red] at (12,-1.6) {$S_{-k}$};
      \path[name path=seg] (5.9,-3.1) -- (7.3,3.1);
      \path[name intersections={of=seg and Smk, by=x}];
      \path[name intersections={of=seg and S0, by=z}];
      \path[name intersections={of=seg and Sk, by=y}];
      \draw[dashed] (5.9,-3.1) -- (x);
      \draw[dashed] (y) -- (7.3,3.1);
      \draw[very thick] (x) -- (y);
      \fill (x) circle (1.3pt) node[below left] {$x_k$};
      \fill (z) circle (1.3pt) node[below left=1pt] {$z_k$};
      \fill (y) circle (1.3pt) node[above left] {$y_k$};
      \coordinate (q) at ($(z)!0.6!(y)$);
      \path[name path=vp] (5.7,-1) -- (5.7,1);
      \path[name intersections={of=vp and S0, by=p}];
      \draw[dotted, thick] (p) -- (q);
      \fill (q) circle (1.1pt) node[above left] {$q$};
      \fill (p) circle (1.1pt) node[above left] {$p$};
      \draw[decorate, decoration={brace, amplitude=4pt, mirror}] ([xshift=6pt]x) -- ([xshift=6pt]z)
        node[midway, right=5pt] {$L_k^-$};
      \draw[decorate, decoration={brace, amplitude=4pt, mirror}] ([xshift=6pt]z) -- ([xshift=6pt]y)
        node[midway, right=5pt] {$L_k^+$};
    \end{tikzpicture}

      (a)
    \end{minipage}%
    \begin{minipage}[b]{0.46\textwidth}
      \centering
    \begin{tikzpicture}[scale=0.42, every node/.style={font=\footnotesize}]
      \draw[thick, name path=S0b] plot[smooth, tension=0.7] coordinates {(0,-0.1) (2.5,0.25) (5,0) (7.5,0.15) (10,0)};
      \draw[thick, gray!70!black, name path=A] plot[smooth, tension=0.7] coordinates {(0,-2.2) (3,-2.0) (6,-2.15) (10,-2.0)};
      \node[right] at (10,0) {$S_0$};
      \node[right, gray!70!black] at (10,-2.0) {$A$};
      \draw[gray!60] plot[smooth] coordinates {(4.15,-3.0) (4.6,0) (5.1,2.2) (5.35,3.6)};
      \draw[gray!60] plot[smooth] coordinates {(4.85,-3.0) (5.3,0) (5.8,2.2) (6.05,3.6)};
      \draw[very thick, name path=gam] plot[smooth] coordinates {(4.5,-3.0) (4.95,0) (5.45,2.2) (5.7,3.6)};
      \path[name intersections={of=gam and S0b, by=zb}];
      \fill (zb) circle (1.3pt) node[below right] {$z$};
      \node[above] at (5.7,3.6) {$\gamma$};
      \node[gray, right] at (6.05,3.6) {$\gamma_k$};
      \coordinate (gt) at (5.45,2.2);
      \fill (gt) circle (1.3pt) node[right=7pt] {$\gamma(t)$};
      \coordinate (a) at (3.2,-2.02);
      \draw[dotted, thick, name path=ac] (a) -- (gt);
      \path[name intersections={of=ac and S0b, by=w}];
      \fill (a) circle (1.1pt);
      \fill (w) circle (1.1pt);
      \draw[decorate, decoration={brace, amplitude=3pt}] ([xshift=-4pt]a) -- ([xshift=-4pt]w)
        node[midway, left=4pt] {$\leq R$};
    \end{tikzpicture}

      (b)
    \end{minipage}
    \caption{(a) The hypersurfaces $S_{-k}$, $S_0$ and $S_k$ need not be ordered by $\preceq$.
      A maximizing segment $\gamma_k$ from $x_k \in P_k$ to $y_k \in F_k$ crosses $S_0$ at $z_k$ (Step~2);
      $L_k^\pm$ are the lengths of its two halves (Step~3), and $p$ and $q$ are as in Step~4.
      \newline (b) The limit curve $\gamma$ through $z$ (Steps~6 and 7)
      and a causal curve from $A$ to $\gamma(t)$ meeting $S_0$ (Step~8).}
    \label{fig:quasi}
  \end{figure}
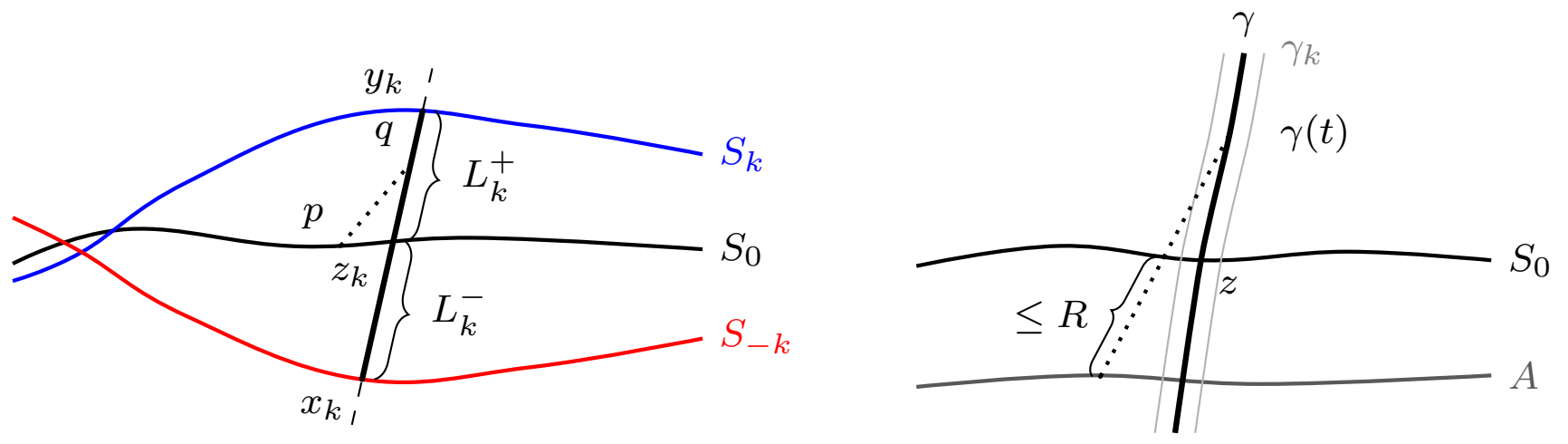

  \emph{Step 3: both halves of $\gamma_k$ have length close to $k$.}
  Let $L_k^- \coloneqq d(x_k, z_k)$ and $L_k^+ \coloneqq d(z_k, y_k)$.
  Then $L_k^- \leq d(P_k, S_0) \leq k + \varepsilon$ and $L_k^+ \leq d(S_0, F_k) \leq k + \varepsilon$.
  Since $\gamma_k$ is maximizing, we have $L_k^- + L_k^+ = L_k \geq 2k - \varepsilon$.
  Hence
  \begin{equation}
    k - 2\varepsilon \leq L_k^\pm \leq k + \varepsilon.
  \end{equation}

  \emph{Step 4: $z_k$ almost realizes the distances from and to $S_0$.}
  Let $q$ be a point of $\gamma_k$ between $z_k$ and $y_k$.
  For $p \in S_0$ with $p \leq q$, we have $d(p, q) + d(q, y_k) \leq d(p, y_k) \leq d(S_0, F_k) \leq k + \varepsilon$.
  For $p \in S_0$ with $p \not\leq q$, we have $d(p, q) = 0 \leq k + \varepsilon - d(q, y_k)$, since $d(q, y_k) \leq L_k^+ \leq k + \varepsilon$.
  Hence $d(S_0, q) \leq k + \varepsilon - d(q, y_k)$.
  Since $d(z_k, q) = L_k^+ - d(q, y_k)$, we obtain
  \begin{equation}\label{eq:defect}
    d(S_0, q) - d(z_k, q) \leq k + \varepsilon - L_k^+ \leq 3\varepsilon.
  \end{equation}
  Similarly, for a point $q$ of $\gamma_k$ between $x_k$ and $z_k$, we have $d(q, S_0) \leq k + \varepsilon - d(x_k, q)$
  and $d(q, z_k) = L_k^- - d(x_k, q)$, so that
  \begin{equation}\label{eq:defect-past}
    d(q, S_0) - d(q, z_k) \leq k + \varepsilon - L_k^- \leq 3\varepsilon.
  \end{equation}

  \emph{Step 5: the endpoints escape.}
  We have $d(S_0, y_k) \geq L_k^+ \geq k - 2\varepsilon \to \infty$ as $k \to \infty$,
  while the function $q \mapsto d(S_0, q)$ is continuous by Lemma~\ref{lem:cauchy} (iii),
  and hence bounded on compact sets.
  Thus $(y_k)_k$ leaves every compact set. Similarly, $(x_k)_k$ leaves every compact set.

  \emph{Step 6: a limit curve.}
  We follow the limit curve argument of~\cite[Lemma 2.2]{Mi}; see also~\cite[p.~298]{EG}.
  Fix a complete Riemannian metric $h$ on $M$.
  Let $\hat\gamma_k$ be the inextendible timelike geodesic extending $\gamma_k$,
  parametrized by $h$-arc length with $\hat\gamma_k(0) = z_k$;
  its domain is $\R$ by~\cite[Lemma 2.17]{MiLR}.
  Then $\gamma_k = \hat\gamma_k|_{[-\alpha_k, \beta_k]}$ with $\hat\gamma_k(-\alpha_k) = x_k$ and $\hat\gamma_k(\beta_k) = y_k$,
  where $\alpha_k$ and $\beta_k$ are at least the $h$-distances from $z_k$ to $x_k$ and to $y_k$, respectively.
  Since $S_0$ is compact, we may assume, by passing to a subsequence, that $z_k \to z \in S_0$.
  Since $x_k$ and $y_k$ leave every compact set by Step~5, we have $\alpha_k, \beta_k \to \infty$.
  By the limit curve theorem~\cite[Theorem 2.51]{MiLR}, applied with the constant sequence of metrics $g$,
  a subsequence of $(\hat\gamma_k)_k$ converges $h$-uniformly on compact subsets of $\R$
  to an inextendible continuous causal curve $\gamma \colon \R \to M$ with $\gamma(0) = z$.
  For every compact interval $[a, b] \subset \R$, we have $[a, b] \subset [-\alpha_k, \beta_k]$ for large $k$,
  so that $\hat\gamma_k|_{[a, b]} = \gamma_k|_{[a, b]}$ is a maximizing segment of $\gamma_k$.

  \emph{Step 7: $\gamma$ is maximizing.}
  For $a < b$, the definition of $d$, the upper semicontinuity of the length~\cite[Theorem 2.41]{MiLR}
  and the continuity of $d$ (Fact~\ref{fact:gh}) give
  \begin{equation}
    \begin{aligned}
      d(\gamma(a), \gamma(b)) & \geq L_g(\gamma|_{[a, b]}) \geq \limsup_{k \to \infty} L_g(\gamma_k|_{[a, b]}) \\
                              & = \limsup_{k \to \infty} d(\gamma_k(a), \gamma_k(b)) = d(\gamma(a), \gamma(b)).
    \end{aligned}
  \end{equation}
  Hence $\gamma$ is an inextendible maximizing continuous causal curve,
  and thus, after reparametrization, an inextendible causal geodesic~\cite[Theorem 2.20]{MiLR},
  which is either timelike or null.
  Fix $t > 0$. For large $k$, the point $\hat\gamma_k(t)$ lies on $\gamma_k$ between $z_k$ and $y_k$,
  so that~\eqref{eq:defect} applies to $q = \hat\gamma_k(t)$.
  Letting $k \to \infty$ and using the continuity of $d(S_0, \cdot)$ and $d$, we obtain the first inequality of
  \begin{equation}\label{eq:defect-limit}
    d(S_0, \gamma(t)) - d(z, \gamma(t)) \leq 3\varepsilon \quad (t > 0), \qquad d(\gamma(t), S_0) - d(\gamma(t), z) \leq 3\varepsilon \quad (t < 0),
  \end{equation}
  and the second one is obtained in the same way from~\eqref{eq:defect-past}.
  Since~\eqref{eq:defect-limit} concerns only the points of $\gamma$ to the future and to the past of $z$,
  it holds in every future-directed parametrization of $\gamma$ with $\gamma(0) = z$.

  \emph{Step 8: $\gamma$ is not null.}
  Choose a smooth acausal compact Cauchy hypersurface $A$ with $A \subset I^-(S_0)$ (see Figure~\ref{fig:quasi} (b)).
  For instance, take a level set $A = \{\mathcal{T} = c\}$ of a Cauchy temporal function $\mathcal{T}$~\cite{BS}, \cite[Theorem 4.120]{MiLR}
  with $c < \min_{S_0} \mathcal{T}$.
  Since $\mathcal{T}$ is a time function, $A$ is disjoint from $S_0 \cup I^+(S_0)$,
  so that $A \subset I^-(S_0)$ by Lemma~\ref{lem:cauchy} (i); moreover, $R \coloneqq d(A, S_0) < \infty$ by Fact~\ref{fact:gh}.
  For $q \in J^+(S_0)$, every causal curve from $A$ to $q$ meets $S_0$ by Lemma~\ref{lem:cauchy} (ii), and hence
  \begin{equation}
    d(A, q) \leq R + d(S_0, q).
  \end{equation}
  By Fact~\ref{fact:levels} applied to $A$, we have $d(A, \gamma(t)) \to \infty$ as $t \to \infty$.
  Suppose that $\gamma$ is null. Then $d(z, \gamma(t)) = 0$ for $t > 0$, since $\gamma$ is maximizing,
  and~\eqref{eq:defect-limit} gives $d(A, \gamma(t)) \leq R + 3\varepsilon$ for all $t > 0$, a contradiction.
  Hence $\gamma$ is a timelike line, and it is complete by timelike geodesic completeness.
  After reparametrizing $\gamma$ by proper time with $\gamma(0) = z$,
  we have $d(z, \gamma(t)) = t$ for $t > 0$ and $d(\gamma(t), z) = -t$ for $t < 0$, since $\gamma$ is maximizing.
  Thus~\eqref{eq:defect-limit} gives the desired estimates.
\end{proof}

\begin{corollary}\label{cor:quasi}
  Suppose that $(M, g)$ satisfies the timelike convergence condition.
  Then $(M, g)$ splits if and only if $(\CM, \dJ)$ contains a $(1, \varepsilon)$-quasi-geodesic line for some $\varepsilon \geq 0$.
\end{corollary}

\begin{proof}
  If $(M, g)$ splits, then the slices $\{t = c\}$, $c \in \R$, form a geodesic line in $(\CM, \dJ)$
  by Proposition~\ref{prop:product} (ii) below.
  Conversely, if $(\CM, \dJ)$ contains a $(1, \varepsilon)$-quasi-geodesic line, then $M$ contains a complete timelike line
  by Theorem~\ref{thm:quasi}, and $(M, g)$ splits by Fact~\ref{fact:splitting}.
\end{proof}

Without the timelike convergence condition, Theorem~\ref{thm:quasi} shows that
a single quasi-geodesic line forces every distance curve to be a quasi-geodesic line.

\begin{corollary}\label{cor:distance-curve}
  Let $(S_t)_{t \in \R}$ be a $(1, \varepsilon)$-quasi-geodesic line in $(\CM, \dJ)$,
  and let $\Sigma$ be an acausal compact Cauchy hypersurface.
  Then $D_\infty(\Sigma) \leq 6\varepsilon + 2\dJ(\Sigma, S_0)$.
  In particular, the distance curve $c_\Sigma$ is a $(1, 6\varepsilon + 2\dJ(\Sigma, S_0))$-quasi-geodesic line.
\end{corollary}

\begin{proof}
  Let $\gamma$ be the complete timelike line of Theorem~\ref{thm:quasi},
  parametrized by proper time with $\gamma(0) = z \in S_0$, and let $R \coloneqq \dJ(\Sigma, S_0)$.
  Since $d_\Sigma$ is a Cauchy time function (Fact~\ref{fact:levels}),
  the function $d_\Sigma \circ \gamma$ is continuous, strictly increasing and takes all real values.
  Moreover, $\abs{d_\Sigma(z)} \leq R$, since $d(\Sigma, z) \leq d(\Sigma, S_0) \leq R$ and $d(z, \Sigma) \leq d(S_0, \Sigma) \leq R$
  by Lemma~\ref{lem:max}.
  Let $k \in \N$ with $k > R$, and let $s_k < 0 < t_k$ be such that $\gamma(s_k) \in \Sigma_{-k}$ and $\gamma(t_k) \in \Sigma_k$.
  By Lemma~\ref{lem:lipschitz} and Theorem~\ref{thm:quasi}, we have
  \begin{equation}
    k = d(\Sigma, \gamma(t_k)) \leq d(S_0, \gamma(t_k)) + R \leq t_k + 3\varepsilon + R,
  \end{equation}
  and similarly $k = d(\gamma(s_k), \Sigma) \leq \abs{s_k} + 3\varepsilon + R$.
  Since $\gamma$ is a timelike line, we obtain
  \begin{equation}
    d(\Sigma_{-k}, \Sigma_k) \geq d(\gamma(s_k), \gamma(t_k)) = t_k - s_k \geq 2k - 6\varepsilon - 2R,
  \end{equation}
  i.e., $\Delta_k \leq 6\varepsilon + 2R$.
  Letting $k \to \infty$, Proposition~\ref{prop:parallel} (ii) gives $D_\infty(\Sigma) \leq 6\varepsilon + 2R$,
  and the last statement follows from Proposition~\ref{prop:parallel} (iii).
\end{proof}

\begin{proposition}\label{prop:lambda}
  Let $\Sigma$ be an acausal compact Cauchy hypersurface.
  Then the rescaled distance curve $t \mapsto \Sigma_{\sqrt{2} t}$ is a $(\sqrt{2}, 0)$-quasi-geodesic line in $(\CM, \dJ)$.
  Consequently, the conclusion of Theorem~\ref{thm:quasi} does not hold for $(\lambda, \varepsilon)$-quasi-geodesic lines
  with $\lambda \geq \sqrt{2}$ in general.
\end{proposition}

\begin{proof}
  Let $\mu > 0$ and $s, t \in \R$. If $s$ and $t$ have the same sign, then
  $\dJ(\Sigma_{\mu s}, \Sigma_{\mu t}) = \mu \abs{t - s}$ by Proposition~\ref{prop:parallel} (i).
  If $s = -\alpha \leq 0 \leq t = \beta$, then Lemma~\ref{lem:lower} gives
  \begin{equation}
    \frac{\mu}{2} \abs{t - s} = \frac{\mu(\alpha + \beta)}{2} \leq \mu \max\{\alpha, \beta\}
    \leq \dJ(\Sigma_{\mu s}, \Sigma_{\mu t}) \leq \mu(\alpha + \beta) = \mu \abs{t - s}.
  \end{equation}
  The case $t \leq 0 \leq s$ follows from the symmetry of $\dJ$.
  Thus $\frac{\mu}{2} \abs{t - s} \leq \dJ(\Sigma_{\mu s}, \Sigma_{\mu t}) \leq \mu \abs{t - s}$ for all $s, t \in \R$.
  For $\mu = \sqrt{2}$, this shows the first claim, since $\sqrt{2}/2 = 1/\sqrt{2}$.
  Applying this to the example of~\cite[Theorem 2]{EG}, which contains no timelike line,
  we obtain $(\lambda, \varepsilon)$-quasi-geodesic lines for all $\lambda \geq \sqrt{2}$ and $\varepsilon \geq 0$
  in a spacetime without complete timelike lines.
\end{proof}

\section{Gromov hyperbolicity and the proof of Theorem~\ref{thm:main}}\label{sec:main}

In this section, we first study $(\CM, \dJ)$ for product spacetimes,
and then show that Gromov hyperbolicity of $(\CM, \dJ)$ forces the defect of every distance curve to be finite.
Combining these with Theorem~\ref{thm:quasi}, we prove Theorem~\ref{thm:main},
and we give examples showing that the timelike convergence condition cannot be dropped.

\subsection{Product spacetimes}

We first consider warped products, which include product spacetimes and the examples of Section~\ref{subsec:examples}.

\begin{lemma}\label{lem:warped}
  Let $(S, h)$ be a compact connected Riemannian manifold, let $a \colon \R \to (0, \infty)$ be a smooth function,
  and let $g = -dt^2 + a(t)^2 h$ on $M = \R \times S$.
  Then $(M, g)$ is globally hyperbolic, each slice $T_c \coloneqq \{c\} \times S$ is a Cauchy hypersurface,
  and the following hold:
  \begin{enumerate}
    \item every Cauchy hypersurface $T$ of $M$ is the graph of a continuous function $f \colon S \to \R$;
    \item for such $T$ and $c \in \R$, we have $\dJ(T, T_c) = \max_S \abs{f - c}$.
  \end{enumerate}
\end{lemma}

\begin{proof}
  Since $S$ is compact, and hence complete, each slice $T_c$ is a Cauchy hypersurface by~\cite[Theorem 3.69 (1)]{BEE};
  in particular, $(M, g)$ is globally hyperbolic.

  (i) Since the slices are Cauchy hypersurfaces and $\partial_t$ is timelike, $t$ is a Cauchy temporal function.
  Hence (i) follows from Step~3 of the proof of Proposition~\ref{prop:proper} and Lemma~\ref{lem:topology}, applied with $\mathcal{T} = t$,
  neither of which uses timelike geodesic completeness.

  (ii) Let $\sigma(s) = (t(s), x(s))$ be a future-directed causal curve from $(t_1, y)$ to $(t_2, x)$.
  Then its length satisfies
  \begin{equation}
    L_g(\sigma) = \int \sqrt{\dot{t}^2 - a(t)^2 \abs{\dot{x}}_h^2} \, ds \leq \int \dot{t} \, ds = t_2 - t_1,
  \end{equation}
  so that $d((t_1, y), (t_2, x)) \leq \max\{t_2 - t_1, 0\}$,
  and the vertical segment shows $d((t_1, x), (t_2, x)) = t_2 - t_1$ if $t_1 \leq t_2$.
  Hence $d(T_c, T) = \max\{\max_S f - c, 0\}$ and $d(T, T_c) = \max\{c - \min_S f, 0\}$,
  and Lemma~\ref{lem:max} gives $\dJ(T, T_c) = \max_S \abs{f - c}$.
\end{proof}

If $(S, h)$ is compact, the product spacetime $(\R \times S, -dt^2 + h)$ is globally hyperbolic and spatially compact by Lemma~\ref{lem:warped} with $a \equiv 1$,
and it is geodesically complete~\cite[Theorem 3.67]{BEE}.

\begin{proposition}\label{prop:product}
  Let $(M, g) = (\R \times S, -dt^2 + h)$, where $(S, h)$ is a compact connected Riemannian manifold,
  and let $D \coloneqq \diam(S, h)$. Then the following hold:
  \begin{enumerate}
    \item every $T \in \CM$ is the graph of a function $f \colon S \to \R$
          satisfying $\abs{f(x) - f(y)} \leq \operatorname{dist}_h(x, y)$ for all $x, y \in S$.
          In particular, $\max f - \min f \leq D$.
          Conversely, the graph of every such function belongs to $\CM$;
    \item for $T = \{t = f(x)\} \in \CM$ and $c \in \R$, we have $\dJ(T, T_c) = \max_S \abs{f - c}$.
          In particular, $\dJ(T_c, T_{c'}) = \abs{c - c'}$ for all $c, c' \in \R$;
    \item the map $\varphi \colon \CM \to \R$, $\varphi(\{t = f(x)\}) \coloneqq (\max f + \min f)/2$, is a $(1, D)$-quasi-isometry;
    \item $(\CM, \dJ)$ is $3D$-hyperbolic;
    \item $(\CM, \dJ)$ is a proper geodesic metric space with at least two ends.
  \end{enumerate}
\end{proposition}

\begin{proof}
  (i) By Lemma~\ref{lem:warped} (i), $T$ is the graph of a function $f \colon S \to \R$.
  The Lipschitz bound follows from the argument of Step~4 of the proof of Proposition~\ref{prop:proper} with $\mathcal{T} = t$ and $\kappa = 1$,
  since $\lambda \partial_t + w$ is future-directed timelike whenever $\lambda > \abs{w}_h$.

  Conversely, let $T$ be the graph of a function $f$ with $\abs{f(x) - f(y)} \leq \operatorname{dist}_h(x, y)$.
  If a future-directed timelike curve $s \mapsto (t(s), \sigma(s))$ joined $(f(x), x)$ to $(f(y), y)$,
  then $\dot t > \abs{\dot\sigma}_h$ would give $f(y) - f(x) = \int \dot t \, ds > \int \abs{\dot\sigma}_h \, ds \geq \operatorname{dist}_h(x, y)$,
  which is impossible. Hence $T$ is achronal.
  Moreover, $t$ takes all real values along every inextendible timelike curve, since the slices are Cauchy hypersurfaces (Lemma~\ref{lem:warped}),
  so that $t - f(x)$ changes its sign along the curve, which therefore meets $T$.

  (ii) By Lemma~\ref{lem:warped} (ii) with $a \equiv 1$, we have $\dJ(T, T_c) = \max_S \abs{f - c}$,
  and applying this to $T = T_{c'}$ gives $\dJ(T_c, T_{c'}) = \abs{c - c'}$.

  (iii) Let $T = \{t = f(x)\}$ and $T' = \{t = f'(x)\}$.
  By (i) and (ii), we have $\dJ(T, T_{\varphi(T)}) = \max_S \abs{f - \varphi(T)} = (\max f - \min f)/2 \leq D/2$,
  and similarly for $T'$.
  The triangle inequality and (ii) give
  \begin{equation}
    \begin{aligned}
      \abs{\dJ(T, T') - \abs{\varphi(T) - \varphi(T')}}
      & = \abs{\dJ(T, T') - \dJ(T_{\varphi(T)}, T_{\varphi(T')})} \\
      & \leq \dJ(T, T_{\varphi(T)}) + \dJ(T', T_{\varphi(T')}) \leq D.
    \end{aligned}
  \end{equation}
  Since $\varphi(T_c) = c$ for every $c \in \R$, the map $\varphi$ is surjective, and hence it is a $(1, D)$-quasi-isometry.

  (iv) By (iii), each of the three distances in the definition of the Gromov product changes by at most $D$ under $\varphi$,
  so that the Gromov products in $\CM$ and in $\R$ satisfy
  $\abs{(x \cdot y)_w - (\varphi(x) \cdot \varphi(y))_{\varphi(w)}} \leq 3D/2$ for all $x, y, w \in \CM$.
  Since $\R$ is $0$-hyperbolic, we obtain, for all $x, y, z, w \in \CM$,
  \begin{equation}
    \begin{aligned}
      (x \cdot y)_w & \geq (\varphi(x) \cdot \varphi(y))_{\varphi(w)} - 3D/2 \\
                    & \geq \min\{(\varphi(x) \cdot \varphi(z))_{\varphi(w)}, (\varphi(z) \cdot \varphi(y))_{\varphi(w)}\} - 3D/2 \\
                    & \geq \min\{(x \cdot z)_w, (z \cdot y)_w\} - 3D.
    \end{aligned}
  \end{equation}

  (v) The space $(\CM, \dJ)$ is a proper geodesic metric space by Theorem~\ref{thm:injective},
  and $\varphi$ is a quasi-isometry to $\R$ by (iii).
  Since $\R$ has two ends, so does $(\CM, \dJ)$~\cite[Chapter I.8, Proposition 8.29]{BH}.
\end{proof}
\subsection{Gromov hyperbolicity and the defect}

\begin{theorem}\label{thm:hyperbolic}
  If $(\CM, \dJ)$ is $\delta$-hyperbolic, then $D_\infty(\Sigma) \leq 8\delta$
  for every acausal compact Cauchy hypersurface $\Sigma$.
  In particular, $M$ contains a complete timelike line.
\end{theorem}

\begin{proof}
  Let $\Sigma$ be an acausal compact Cauchy hypersurface; such hypersurfaces exist by~\cite{BS}, \cite[Theorem 4.120]{MiLR}
  (see Section~\ref{sec:prelim}).
  We put $a_t \coloneqq \Sigma_t$ and $b_t \coloneqq \Sigma_{-t}$ for $t \geq 0$.
  Fix $r > 4\delta$ and let $s \geq r$.
  Applying the four-point condition of Definition~\ref{def:hyperbolic} twice, we have
  \begin{equation}
    \begin{aligned}
      (a_r \cdot b_r)_\Sigma & \geq \min\{(a_r \cdot a_s)_\Sigma, (a_s \cdot b_r)_\Sigma\} - \delta \\
                     & \geq \min\{(a_r \cdot a_s)_\Sigma, (a_s \cdot b_s)_\Sigma, (b_s \cdot b_r)_\Sigma\} - 2\delta
                       = \min\{r, (a_s \cdot b_s)_\Sigma\} - 2\delta,
    \end{aligned}
  \end{equation}
  where the last equality follows from Lemma~\ref{lem:gromov}.
  By Lemma~\ref{lem:gromov} again, the left-hand side is at most $r/2$.
  If $(a_s \cdot b_s)_\Sigma \geq r$, then we would have $r/2 \geq r - 2\delta$, i.e., $r \leq 4\delta$, a contradiction.
  Hence $(a_s \cdot b_s)_\Sigma < r$, and we obtain $(a_s \cdot b_s)_\Sigma \leq r/2 + 2\delta$.
  By Lemma~\ref{lem:gromov}, this means $D(s, s) \leq r + 4\delta$ for all $s \geq r$.
  By Proposition~\ref{prop:parallel} (ii), we obtain $D_\infty(\Sigma) = \lim_{k \to \infty} D(k, k) \leq r + 4\delta$.
  Letting $r \to 4\delta$, we conclude $D_\infty(\Sigma) \leq 8\delta$.
  The last statement follows from Proposition~\ref{prop:parallel} (iii) and Theorem~\ref{thm:quasi}.
\end{proof}

The following fact is well known; we include a proof for completeness.

\begin{lemma}\label{lem:ends}
  A proper geodesic metric space $X$ with at least two ends contains a geodesic line.
\end{lemma}

\begin{proof}
  By~\cite[Chapter I.8, Lemma 8.28 (2)]{BH}, there are geodesic rays $r_1$ and $r_2$ which do not converge to the same end.
  Then there exists a compact set $K \subset X$ such that, for every $N \geq 0$,
  the sets $r_1([N, \infty))$ and $r_2([N, \infty))$ do not lie in the same path component of $X \setminus K$.
  For $i = 1, 2$, we have $\rho(r_i(t), K) \geq t - \max_{x \in K} \rho(r_i(0), x) \to \infty$ as $t \to \infty$.
  Hence there exists $N_0$ such that $r_1([N_0, \infty))$ and $r_2([N_0, \infty))$ are contained in $X \setminus K$.
  Since these sets are path connected, they lie in path components $U_1$ and $U_2$ of $X \setminus K$, respectively,
  and $U_1 \neq U_2$.
  Let $p_n \coloneqq r_1(N_0 + n)$ and $q_n \coloneqq r_2(N_0 + n)$, so that $\rho(p_n, K) \to \infty$ and $\rho(q_n, K) \to \infty$.
  A geodesic segment from $p_n$ to $q_n$ meets $K$ at a point $o_n$,
  since otherwise it would be a path in $X \setminus K$ from $p_n \in U_1$ to $q_n \in U_2$,
  and $U_1$ and $U_2$ would be the same path component.
  We parametrize it by arc length with $o_n$ at time $0$;
  then its domain contains $[-\rho(p_n, K), \rho(q_n, K)]$.
  Since $o_n$ lies in the compact set $K$,
  the Arzel\`a--Ascoli theorem~\cite[Chapter I.3, Lemma 3.10]{BH} and a diagonal subsequence argument give a subsequence
  (as in the proof of~\cite[Chapter I.8, Lemma 8.28 (2)]{BH})
  converging uniformly on compact sets to a curve $\R \to X$ which is isometric on every compact interval,
  i.e., a geodesic line.
\end{proof}

\begin{remark}[Clock defects]\label{rem:clock}
  The defect $D_\infty(\Sigma)$ has a physical meaning~\cite{Mi}.
  For a complete future-directed timelike geodesic $\gamma$ parametrized by proper time,
  the function $t \mapsto d_\Sigma(\gamma(t)) - t$ is non-decreasing~\cite[Lemma 7.1]{Mi},
  so that the limits $A^\pm_\Sigma(\gamma) \coloneqq \lim_{t \to \pm\infty} (d_\Sigma(\gamma(t)) - t)$ exist~\cite[Corollary 7.5]{Mi}.
  We call $\cdf_\Sigma(\gamma) \coloneqq A^+_\Sigma(\gamma) - A^-_\Sigma(\gamma) \in [0, \infty]$
  the \emph{total clock defect} of $\gamma$, following the interpretation of these limits as clock defects
  in~\cite[Corollary 7.5 and the discussion after Theorem 4.5]{Mi};
  it compares the proper time of the observer $\gamma$ with the signed distance from $\Sigma$ over the whole history of $\gamma$.
  By~\cite[Corollary 7.5]{Mi}, we have $A^+_\Sigma(\gamma) \geq b^+(\gamma(0))$ and $A^-_\Sigma(\gamma) \leq b^-(\gamma(0))$,
  so that $\cdf_\Sigma(\gamma) \geq (b^+ - b^-)(\gamma(0)) \geq D_\infty(\Sigma)$ by Proposition~\ref{prop:minguzzi}.
  If $D_\infty(\Sigma) < \infty$, then $\dom b^+ \cap \dom b^- \neq \emptyset$ by Proposition~\ref{prop:minguzzi},
  and through a minimum point of $b^+ - b^-$ there passes a complete timelike line $\gamma$
  which is calibrated for both $b^+$ and $b^-$~\cite[Theorem 7.2 and its proof]{Mi};
  for this $\gamma$, both inequalities are equalities~\cite[Corollary 7.5]{Mi}.
  Hence $D_\infty(\Sigma) = \inf_\gamma \cdf_\Sigma(\gamma)$, the infimum being attained if $D_\infty(\Sigma) < \infty$,
  and, without the timelike convergence condition,
  Theorem~\ref{thm:hyperbolic} and Corollary~\ref{cor:distance-curve} bound the clock defect of some observer
  by constants of coarse geometry.
\end{remark}

\subsection{Proof of Theorem~\ref{thm:main}}

\begin{theorem}\label{thm:main}
  Let $(M, g)$ be a spatially compact, globally hyperbolic spacetime.
  Suppose that $(M, g)$ is timelike geodesically complete and satisfies the timelike convergence condition.
  Then the following conditions are equivalent:
  \begin{enumerate}
    \item $(M, g)$ splits isometrically as $(\R \times S, -dt^2 + h)$,
          where $(S, h)$ is a compact Riemannian manifold with non-negative Ricci curvature;
    \item $(\CM, \dJ)$ contains a $(1, \varepsilon)$-quasi-geodesic line for some $\varepsilon \geq 0$;
    \item $(\CM, \dJ)$ is quasi-isometric to $\R$;
    \item $(\CM, \dJ)$ is Gromov hyperbolic;
    \item $(\CM, \dJ)$ has at least two ends.
  \end{enumerate}
\end{theorem}

\begin{proof}
  The equivalence (i) $\Leftrightarrow$ (ii) is Corollary~\ref{cor:quasi}.
  If (i) holds, then (iii), (iv) and (v) hold by Proposition~\ref{prop:product} (iii), (iv) and (v).
  Suppose that (iii) holds.
  For geodesic spaces, Gromov hyperbolicity in the sense of Definition~\ref{def:hyperbolic}
  agrees with the hyperbolicity defined by slim geodesic triangles~\cite[Chapter III.H, Propositions 1.17 and 1.22]{BH}.
  Since $(\CM, \dJ)$ is a geodesic space by Theorem~\ref{thm:injective} and $\R$ is $0$-hyperbolic,
  $(\CM, \dJ)$ is Gromov hyperbolic by the quasi-isometry invariance of hyperbolicity~\cite[Chapter III.H, Theorem 1.9]{BH};
  thus (iv) holds.
  If (iv) holds, then $D_\infty(\Sigma) < \infty$ for every acausal compact Cauchy hypersurface $\Sigma$
  by Theorem~\ref{thm:hyperbolic}, so that the distance curve $c_\Sigma$ is a $(1, D_\infty(\Sigma))$-quasi-geodesic line
  by Proposition~\ref{prop:parallel} (iii); thus (ii) holds.
  If (v) holds, then $(\CM, \dJ)$ is a proper geodesic metric space by Theorem~\ref{thm:injective},
  and it contains a geodesic line by Lemma~\ref{lem:ends}; thus (ii) holds.
\end{proof}

We equip $\CM$ with the topology of uniform convergence of graph functions, which is the topology of $\dJ$ by Lemma~\ref{lem:topology}.
Let $Q \coloneqq [0, 1]^{\N}$ be the Hilbert cube, i.e., the product of countably many copies of $[0, 1]$ with the product topology,
which is a compact metrizable space,
and for $p \in M$, let $\CM(p) \coloneqq \{\Sigma \in \CM ; p \in \Sigma\}$.
Recall that $n \geq 1$, i.e., the spatial dimension is positive.

\begin{corollary}\label{cor:type}
  Under the assumptions of Theorem~\ref{thm:main}, the following conditions are equivalent:
  \begin{enumerate}
    \item $(M, g)$ splits isometrically as $(\R \times S, -dt^2 + h)$;
    \item $\CM$ is homeomorphic to $\R \times Q$;
    \item $\CM(p)$ is compact for some $p \in M$;
    \item $\CM(p)$ is homeomorphic to $Q$ for every $p \in M$.
  \end{enumerate}
  Moreover, if $(M, g) = (\R \times S, -dt^2 + h)$ and $p = (t_0, x_0)$,
  then there exists a homeomorphism of pairs $(\CM, \CM(p)) \cong (\R \times Q, \{t_0\} \times Q)$.
\end{corollary}

\begin{proof}
  We use the following observation.
  Let $X$ be a topological space, and suppose that there exist a continuous function $e \colon X \to \R$ such that $e^{-1}(0)$ is compact
  and proper rays $r_\pm$ in $X$ such that $\pm e(r_\pm(t)) > 0$ for large $t$.
  Then $X$ has at least two ends.
  Indeed, let $N \geq 0$. The sets $r_+([N, \infty))$ and $r_-([N, \infty))$ contain points $r_+(t)$ with $e > 0$ and $r_-(t')$ with $e < 0$,
  and by the intermediate value theorem no path in $X \setminus e^{-1}(0)$ joins two such points.
  Hence, for every $N \geq 0$, the sets $r_+([N, \infty))$ and $r_-([N, \infty))$ do not lie in the same path component of $X \setminus e^{-1}(0)$,
  so that $r_+$ and $r_-$ do not converge to the same end.

  (i) $\Rightarrow$ (ii), (iv) and the last statement.
  Let $(M, g) = (\R \times S, -dt^2 + h)$, $D \coloneqq \diam(S, h)$ and $p = (t_0, x_0)$, and let
  \begin{equation}
    K_{x_0} \coloneqq \{u \colon S \to \R ; \abs{u(x) - u(y)} \leq \operatorname{dist}_h(x, y) \text{ for all } x, y \in S, \ u(x_0) = 0\}
  \end{equation}
  with the uniform topology.
  By Proposition~\ref{prop:product} (i) and Lemma~\ref{lem:topology}, applied with $\mathcal{T} = t$,
  the map $\{t = f(x)\} \mapsto (f(x_0), f - f(x_0))$ is a homeomorphism $\CM \to \R \times K_{x_0}$,
  with inverse $(c, u) \mapsto \{t = c + u(x)\}$, and it maps $\CM(p)$ onto $\{t_0\} \times K_{x_0}$.
  We show that $K_{x_0}$ is homeomorphic to $Q$.
  To this end, we show that $K_{x_0}$ is a convex, compact and infinite-dimensional subset of $C^0(S)$, and then apply Keller's theorem.

  The convexity follows, since the conditions defining $K_{x_0}$ are preserved under convex combinations.

  The compactness follows from the Arzel\`a--Ascoli theorem.
  Indeed, $K_{x_0}$ is closed in $C^0(S)$, since the conditions defining it pass to uniform limits.
  Moreover, $K_{x_0}$ is equicontinuous, since its elements are $1$-Lipschitz,
  and uniformly bounded, since $\abs{u(x)} = \abs{u(x) - u(x_0)} \leq \operatorname{dist}_h(x, x_0) \leq D$ for $u \in K_{x_0}$ and $x \in S$.

  The infinite-dimensionality follows, since $K_{x_0}$ contains an $m$-dimensional simplex for every $m \in \N$.
  Indeed, since $\dim S = n \geq 1$, there exist distinct points $y_1, \dots, y_m \in S \setminus \{x_0\}$ and $r > 0$
  such that the balls $B(y_k, r)$ with respect to $\operatorname{dist}_h$ are pairwise disjoint and do not contain $x_0$.
  The functions $\varphi_k \coloneqq \max\{0, r - \operatorname{dist}_h(\cdot, y_k)\}$ are $1$-Lipschitz and vanish at $x_0$,
  so that they belong to $K_{x_0}$.
  They are linearly independent, since $\varphi_k(y_k) = r > 0$ and $\varphi_k$ vanishes outside $B(y_k, r)$.
  Since $K_{x_0}$ is convex and contains $0$, it contains the $m$-dimensional simplex spanned by $0, \varphi_1, \dots, \varphi_m$.

  Since $C^0(S)$ is a Banach space, $K_{x_0}$ is an infinite-dimensional compact convex subset of a Fr\'echet space.
  By Keller's theorem~\cite{Ke31} in the form of~\cite[Chapter III, Proposition 3.1 and Theorem 3.1]{BP75},
  $K_{x_0}$ is homeomorphic to $Q$.

  Let $\kappa \colon K_{x_0} \to Q$ be a homeomorphism.
  Composing the homeomorphism $\CM \to \R \times K_{x_0}$ above with $\mathrm{id}_\R \times \kappa$,
  we obtain a homeomorphism $\CM \to \R \times Q$ which maps $\CM(p)$ onto $\{t_0\} \times Q$.
  This proves (ii) and the last statement, and (iv) follows, since $p$ is arbitrary.

  (iv) $\Rightarrow$ (iii). This is clear.

  (ii) $\Rightarrow$ (i).
  The observation applies to $X = \R \times Q$ with $e(s, q) \coloneqq s$ and $r_\pm(t) \coloneqq (\pm t, q_0)$ for some $q_0 \in Q$.
  Indeed, $e^{-1}(0) = \{0\} \times Q$ is compact, and $r_\pm$ are proper rays,
  since every compact subset $C$ of $\R \times Q$ is contained in $[-R, R] \times Q$ for some $R > 0$,
  so that $r_\pm^{-1}(C)$ is a closed subset of $[0, R]$.
  Since ends are preserved by homeomorphisms,
  $(\CM, \dJ)$ has at least two ends, and (i) follows from Theorem~\ref{thm:main}.

  (iii) $\Rightarrow$ (i).
  With $\mathcal{T}$ and $H_t$ as before Lemma~\ref{lem:topology}, write $p = (t_0, x_0)$,
  and let $e_p(S) \coloneqq f_S(x_0) - t_0$ for $S \in \CM$.
  Then $e_p$ is continuous by Lemma~\ref{lem:topology}, $e_p^{-1}(0) = \CM(p)$ is compact, and $e_p(H_{\pm t}) = \pm t - t_0$.
  The maps $r_\pm(t) \coloneqq H_{\pm t}$, $t \geq 0$, are continuous by Lemma~\ref{lem:topology}, since $f_{H_t} \equiv t$.
  They are proper rays, since $\max_{S'} \abs{f_S} \leq R$ on every compact subset $C$ of $\CM$ for some $R > 0$ by Lemma~\ref{lem:topology},
  so that $r_\pm^{-1}(C)$ is a closed subset of $[0, R]$.
  By the observation, $(\CM, \dJ)$ has at least two ends, and (i) follows from Theorem~\ref{thm:main}.
\end{proof}

\begin{remark}\label{rem:point}
  We use the notation of the proof of Corollary~\ref{cor:type},
  and for $p \in M$, let $N(p) \coloneqq M \setminus (I^+(p) \cup I^-(p))$ be the set of points which are not chronologically related to $p$. 
  For every $p \in M$, we have $\bigcup_{\Sigma \in \CM(p)} \Sigma = N(p)$,
  and $\CM(p)$ is compact if and only if $N(p)$ is compact.
  Indeed, $\subset$ follows from achronality.
  Conversely, for $q \in N(p)$, choose $c$ with $p, q \in I^+(H_c)$;
  then $A \coloneqq I^-(p) \cup I^-(q) \cup I^-(H_c)$ satisfies the assumptions of Lemma~\ref{lem:lattice} (ii),
  and $p, q \in \partial A \in \CM$, which gives $\supset$.

  If $N(p)$ is compact, then $\CM(p) = e_p^{-1}(0)$ is closed and $\dJ$-bounded by $\max_{N(p)} \abs{d_{\Sigma}}$
  around a fixed acausal compact Cauchy hypersurface $\Sigma$, and hence compact by Proposition~\ref{prop:proper}.
  If $\CM(p)$ is compact, then it is compact with respect to $d_H$ by Lemma~\ref{lem:topology}.
  Given $x_k \in T_k \in \CM(p)$, a subsequence of $(T_k)$ converges to some $T \in \CM(p)$ with respect to $d_H$.
  Since $T$ is compact, there exist $y_k \in T$ with $\rho_0(x_k, y_k) = \rho_0(x_k, T) \leq d_H(T_k, T)$,
  and we may assume $y_k \to y \in T$.
  Then $\rho_0(x_k, y) \leq \rho_0(x_k, y_k) + \rho_0(y_k, y) \to 0$, so that $x_k \to y \in T \subset N(p)$.
  Hence $N(p)$ is compact.
  Thus condition (iii) is equivalent to the point condition of Bartnik~\cite{Ba}, i.e., the compactness of $N(p)$ for some $p$.
  Under this condition, Bartnik proved the existence of a constant mean curvature Cauchy hypersurface through $p$~\cite[Theorem]{Ba},
  and hence the splitting under the assumptions of Theorem~\ref{thm:main}~\cite[Corollary 2]{Ba}.
\end{remark}

The number of ends of $\CM$ is at most two, without any curvature or completeness condition.

\begin{proposition}\label{prop:two-ends}
  Let $(M, g)$ be a spatially compact, globally hyperbolic spacetime,
  and let $\mathcal{T}$ and $H_t$ be as before Lemma~\ref{lem:topology},
  where $\mathcal{T}$ may be any smooth Cauchy temporal function as in Step~3 of the proof of Proposition~\ref{prop:proper}.
  Then $\CM$ has either one or two ends.
  Moreover, it has two ends if and only if there exist $a < b$ such that no Cauchy hypersurface meets both $I^-(H_a)$ and $I^+(H_b)$.
\end{proposition}

\begin{proof}
  By Lemma~\ref{lem:lattice} (i) and the remark after Lemma~\ref{lem:topology},
  $S \preceq T$ if and only if $f_S \leq f_{T}$,
  and $S \vee T$ and $S \wedge T$ are the graphs of $\max\{f_S, f_{T}\}$ and $\min\{f_S, f_{T}\}$.

  \emph{Step 1: order intervals.}
  For $a < b$, the order interval $[H_a, H_b] = \{S ; a \leq f_S \leq b\}$ is compact.
  Indeed, Steps~4--7 of the proof of Proposition~\ref{prop:proper} use only $S \subset [a, b] \times S'$, and hence apply to sequences in $[H_a, H_b]$.
  Thus its elements are graphs of uniformly Lipschitz functions (Step~4), a subsequence of their graph functions converges uniformly (Step~5),
  the limit graph belongs to $\CM$ (Step~6), and hence to $[H_a, H_b]$, and the subsequence converges to it (Step~7).
  Conversely, every compact set $K \subset \CM$ is contained in $[H_a, H_b]$ with $a \coloneqq \min_{S \in K} \min f_S$ and $b \coloneqq \max_{S \in K} \max f_S$,
  which exist since $S \mapsto \min f_S$ and $S \mapsto \max f_S$ are continuous by Lemma~\ref{lem:topology}.
  Hence $t \mapsto H_t$, $t \geq 0$, is a proper ray, so that $\CM$ has at least one end.
  Moreover, every proper ray eventually leaves every $[H_a, H_b]$, since $[H_a, H_b]$ is compact.

  \emph{Step 2: the complement of an order interval.}
  Let $F_b \coloneqq \{S ; \max f_S > b\}$ and $P_a \coloneqq \{S ; \min f_S < a\}$,
  i.e., the sets of Cauchy hypersurfaces meeting $I^+(H_b)$ and $I^-(H_a)$, respectively.
  They are open, and $\CM \setminus [H_a, H_b] = F_b \cup P_a$.
  The set $F_b$ is path connected, since for $S \in F_b$ the path $t \mapsto S \vee H_t$ is continuous, stays in $F_b$,
  and equals $H_t$ for $t \geq \max f_S$, and the path $t \mapsto H_t$, $t > b$, lies in $F_b$.
  Similarly, $P_a$ is path connected.
  Hence $\CM \setminus [H_a, H_b]$ is path connected if $F_b \cap P_a \neq \emptyset$,
  and it is the disjoint union of the open path connected sets $F_b$ and $P_a$ otherwise.

  \emph{Step 3: conclusion.}
  For $n \in \N$, let $K_n \coloneqq [H_{-n}, H_n]$.
  By Step~1, every compact set is contained in some $K_n$, and every proper ray eventually leaves every $K_n$.
  Since $F_n \subset F_b$ and $P_{-n} \subset P_a$ for $n \geq \max\{-a, b\}$,
  the condition $F_b \cap P_a = \emptyset$ holds for some $a < b$ if and only if $F_n \cap P_{-n} = \emptyset$ for some $n$.

  Suppose first that $F_n \cap P_{-n} \neq \emptyset$ for all $n$.
  Let $r_1$ and $r_2$ be proper rays, and let $K$ be compact with $K \subset K_n$.
  Tails of $r_1$ and $r_2$ lie in $\CM \setminus K_n$, which is path connected by Step~2 and contained in $\CM \setminus K$.
  Hence $r_1$ and $r_2$ converge to the same end,
  and $\CM$ has exactly one end, since $t \mapsto H_t$, $t \geq 0$, is a proper ray by Step~1.

  Suppose next that $F_m \cap P_{-m} = \emptyset$ for some $m$.
  Then $F_n \cap P_{-n} = \emptyset$ for all $n \geq m$, since $F_n \subset F_m$ and $P_{-n} \subset P_{-m}$.
  Let $r$ be a proper ray. A tail of $r$ lies in $\CM \setminus K_m = F_m \cup P_{-m}$,
  and since it is connected, it lies in $F_m$ or in $P_{-m}$;
  we say that $r$ is of \emph{future type} or of \emph{past type}, respectively.
  The type persists when $K_m$ is enlarged, which allows us to separate tails from arbitrary compact sets.
  Indeed, if $r$ is of future type and $n \geq m$, then a tail of $r$ lies in $F_n \cup P_{-n}$ by Step~1,
  and hence in $F_n$, since $P_{-n} \subset P_{-m}$ is disjoint from $F_m$;
  similarly, tails of rays of past type lie in $P_{-n}$.

  Let $r_1$ and $r_2$ be proper rays.
  If they are of the same type, say of future type, and $K \subset K_n$ is compact with $n \geq m$,
  then tails of $r_1$ and $r_2$ lie in $F_n$, which is path connected by Step~2 and contained in $\CM \setminus K$,
  so that $r_1$ and $r_2$ converge to the same end.
  If $r_1$ is of future type and $r_2$ is of past type, then for every $N \geq 0$,
  the sets $r_1([N, \infty))$ and $r_2([N, \infty))$ meet $F_m$ and $P_{-m}$, respectively,
  whereas every path in $\CM \setminus K_m$ stays in one of the disjoint open sets $F_m$ and $P_{-m}$.
  Hence $r_1$ and $r_2$ do not converge to the same end.
  Since $t \mapsto H_t$ and $t \mapsto H_{-t}$, $t \geq 0$, are proper rays of future type and of past type, respectively,
  $\CM$ has exactly two ends.
  Since $F_b \cap P_a = \emptyset$ means that no Cauchy hypersurface meets both $I^-(H_a)$ and $I^+(H_b)$,
  this completes the proof.
\end{proof}

We compare the two-end condition with causal conditions from the literature.
Following~\cite{GKP}, a \emph{TIP} (resp.\ \emph{TIF}) is the chronological past $I^-(\gamma)$ (resp.\ future $I^+(\gamma)$)
of a future (resp.\ past) inextendible timelike curve $\gamma$,
and it is \emph{minimal} if it contains no other TIP (resp.\ TIF) as a proper subset.
The \emph{min-min condition} of~\cite{GaHe25} requires that some minimal TIP meets some minimal TIF;
under spatial compactness, minimal TIPs and TIFs exist~\cite[Theorem 1.2]{GaHe25},
and the min-min condition implies the existence of a timelike line~\cite[Theorem 1.4]{GaHe25}.

\begin{proposition}\label{prop:tip}
  Let $(M, g)$ be a spatially compact, globally hyperbolic spacetime.
  If $\CM$ has two ends, then every TIP meets every TIF.
  In particular, the min-min condition holds.
\end{proposition}

\begin{proof}
  We use the notation of Proposition~\ref{prop:two-ends}, and choose $a < b$ such that no Cauchy hypersurface meets both $I^-(H_a)$ and $I^+(H_b)$.
  Let $p \in I^-(H_a)$ and $q \in I^+(H_b)$, so that $\mathcal{T}(p) < a < b < \mathcal{T}(q)$.
  We claim that $p \ll q$. Since $\mathcal{T}$ increases along future-directed timelike curves, $q \not\ll p$.
  Suppose that $p \not\ll q$, choose $c < \mathcal{T}(p)$, and let $A \coloneqq I^-(p) \cup I^-(q) \cup I^-(H_c)$.
  As in Remark~\ref{rem:point}, $S \coloneqq \partial A \in \CM$ by Lemma~\ref{lem:lattice} (ii),
  and $p, q \in S$, since $p, q \notin A$.
  Then $S$ meets both $I^-(H_a)$ and $I^+(H_b)$, a contradiction.
  Now let $I^-(\gamma)$ be a TIP and $I^+(\eta)$ a TIF.
  Since $\mathcal{T}$ is a Cauchy temporal function, there are $q$ on $\gamma$ with $\mathcal{T}(q) > b$ and $p$ on $\eta$ with $\mathcal{T}(p) < a$.
  Then $p \ll q$, and $\emptyset \neq I^+(p) \cap I^-(q) \subset I^+(\eta) \cap I^-(\gamma)$.
\end{proof}

\begin{figure}[t]
  \centering
  \begin{tikzpicture}[scale=1.1]
    \fill[gray!30] (0,1) -- (-2,-1) -- (2,-1) -- cycle;
    \fill[gray!12] (-2,1) -- (0,1) -- (-2,-1) -- cycle;
    \fill[gray!12] (2,1) -- (0,1) -- (2,-1) -- cycle;
    \draw (-2,1) -- (2,1);
    \draw (-2,-1) -- (2,-1);
    \draw[dashed] (-2,-1) -- (-2,1);
    \draw[dashed] (2,-1) -- (2,1);
    \draw[thick] (0,1) -- (-2,-1);
    \draw[thick] (0,1) -- (2,-1);
    \draw[thick] (0.25,-0.85) .. controls (0.33,-0.2) and (-0.12,0.3) .. (0,1);
    \draw[thick] (2,-1) .. controls (1.75,-0.4) and (1.82,0.2) .. (1.6,0.85);
    \filldraw[fill=white] (0,1) circle (2pt) node[above] {$p$};
    \filldraw[fill=white] (2,-1) circle (2pt) node[below right] {$q$};
    \filldraw[fill=white] (-2,-1) circle (2pt) node[below left] {$q$};
    \node at (0.5,-0.75) {$\gamma_1$};
    \node at (1.55,-0.1) {$\gamma_2$};
    \node at (-0.75,-0.6) {$I^-(\gamma_1)$};
    \node at (-1.45,0.55) {$I^+(\gamma_2)$};
    \node[left] at (-2.35,1) {$\pi/2$};
    \node[left] at (-2.35,-1) {$-\pi/2$};
    \node[left] at (-2.35,0) {$\eta$};
    \node[below] at (0,-1.05) {$\theta = 0$};
  \end{tikzpicture}
  \caption{The two-dimensional de Sitter spacetime in the conformal coordinates $(\eta, \theta)$,
    where the dashed sides $\theta = \pm \pi$ are identified.
    The TIP of a curve $\gamma_1$ approaching $p = (\pi/2, 0)$ in the future
    and the TIF of a curve $\gamma_2$ approaching the antipodal point $q = (-\pi/2, \pi)$ in the past
    are separated by the two null geodesics from $q$ to $p$.}
  \label{fig:desitter}
\end{figure}
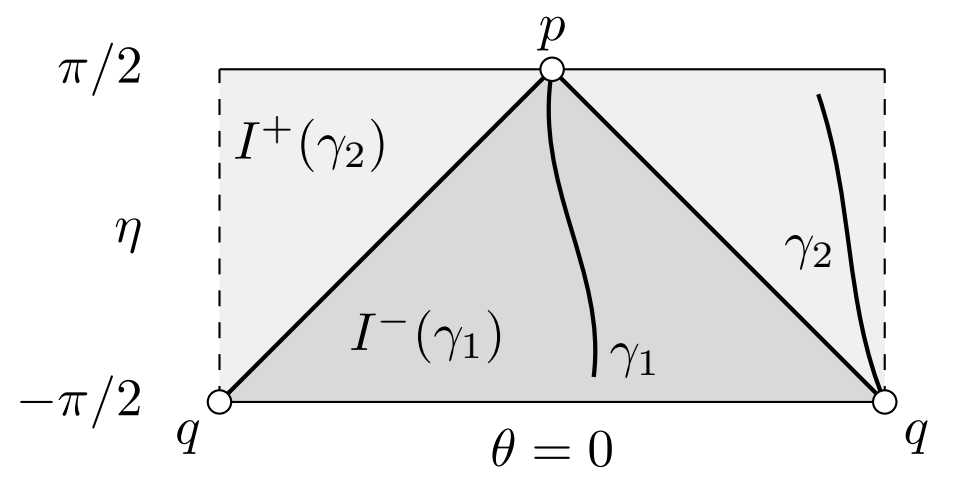

\begin{remark}\label{rem:causal}
  Assume our standing assumptions.
  If $I^-(\gamma) = M$ for every future inextendible causal curve $\gamma$,
  i.e., the no future observer horizons condition~\cite[Definition 2.6]{GHZ26} holds,
  then $N(p)$ is compact for some $p \in M$ by~\cite{Ti93} (see also~\cite[Proposition 3]{GL18}).
  Under this condition, the group of time-orientation preserving isometries acts properly on $M$
  and splits as a semi-direct product of a compact Lie group and a group which is trivial, $\Z$ or $\R$~\cite[Theorem 1.1 and Corollary 1.3]{GHZ26}.
  By Remark~\ref{rem:point} and the proof of (iii) $\Rightarrow$ (i) of Corollary~\ref{cor:type}, which does not use the timelike convergence condition,
  the compactness of $N(p)$ implies that $\CM$ has two ends.
  Together with Proposition~\ref{prop:tip}, we obtain
  \begin{equation}
    \begin{gathered}
      I^-(\gamma) = M \text{ for all } \gamma \Rightarrow N(p) \text{ compact for some } p \Rightarrow \CM \text{ has two ends} \\
      \Rightarrow \text{every TIP meets every TIF} \Rightarrow \text{min-min}.
    \end{gathered}
  \end{equation}
  Under the assumptions of Theorem~\ref{thm:main}, all of these conditions are equivalent to the splitting.
  Without the timelike convergence condition, the second implication is strict by Example~\ref{ex:cylinder}.
  Moreover, the min-min condition does not imply that $\CM$ has two ends,
  as the de Sitter spacetime of Example~\ref{ex:desitter} shows.
  Via $\eta \coloneqq \arctan(\sinh \tau)$, it is conformal to $(-\pi/2, \pi/2) \times \mathbf{S}^n$ in the Einstein static universe $(\R \times \mathbf{S}^n, -d\eta^2 + h_{\mathbf{S}^n})$~\cite[Section 5.2]{HE}.
  As is well known, every TIP is the chronological past of an inextendible timelike curve
  approaching a point $(\pi/2, u_0)$ of the future conformal boundary, and it depends only on $u_0$;
  dually, every TIF is determined by a point $(-\pi/2, u_0)$ of the past conformal boundary,
  and all of them are minimal.
  Taking the same point $u_0$ gives the min-min condition,
  while antipodal points give a TIP and a TIF which do not meet (see Figure~\ref{fig:desitter}),
  so that $\CM$ has one end by Propositions~\ref{prop:tip} and~\ref{prop:two-ends}.
\end{remark}

Under our standing assumptions, the number of ends is determined as follows.

\begin{corollary}\label{cor:ends}
  Let $\Sigma$ be an acausal compact Cauchy hypersurface. Then the following hold:
  \begin{enumerate}
    \item if $\dom b^+ \cap \dom b^- = \emptyset$, then $(\CM, \dJ)$ has exactly one end;
    \item if $(M, g)$ satisfies the timelike convergence condition,
          then $(\CM, \dJ)$ has exactly two ends if $(M, g)$ splits, and exactly one end otherwise.
  \end{enumerate}
\end{corollary}

\begin{proof}
  (i) Suppose that $(\CM, \dJ)$ has at least two ends.
  Then it contains a geodesic line by Theorem~\ref{thm:injective} and Lemma~\ref{lem:ends},
  so that $D_\infty(\Sigma) < \infty$ by Corollary~\ref{cor:distance-curve},
  and $\dom b^+ \cap \dom b^- \neq \emptyset$ by Proposition~\ref{prop:minguzzi} and Fact~\ref{fact:defect},
  contrary to the assumption.
  Hence $(\CM, \dJ)$ has at most one end, and exactly one end by Proposition~\ref{prop:two-ends}.
  (ii) If $(M, g)$ splits, then $(\CM, \dJ)$ has at least two ends by Proposition~\ref{prop:product} (v),
  and hence exactly two ends by Proposition~\ref{prop:two-ends}.
  If $(M, g)$ does not split, then $\dom b^+ \cap \dom b^- = \emptyset$ by Fact~\ref{fact:defect}, and (i) applies.
\end{proof}

\begin{example}\label{ex:EG}
  The example of Ehrlich and Galloway~\cite{EG} is globally hyperbolic, causally geodesically complete
  and spatially compact, and contains no timelike line.
  Hence (ii) of Theorem~\ref{thm:main} fails by Theorem~\ref{thm:quasi}.
  Since the implications (iv) $\Rightarrow$ (ii) and (v) $\Rightarrow$ (ii) in the proof of Theorem~\ref{thm:main}
  do not use the timelike convergence condition, $(\CM, \dJ)$ is not Gromov hyperbolic and has at most one end.
  Since a distance ray $c_\Sigma|_{[0, \infty)}$ is a geodesic ray (Proposition~\ref{prop:parallel} (i)),
  $(\CM, \dJ)$ has exactly one end, so that the future and past distance rays go to the same end.
\end{example}

\subsection{Examples without the timelike convergence condition}\label{subsec:examples}

In this subsection, we show that the timelike convergence condition cannot be dropped in Theorem~\ref{thm:main}.
Example~\ref{ex:warped} satisfies (iii) and (iv) but not (i),
and Example~\ref{ex:desitter} satisfies (ii) but neither (iii) nor (iv).
Both examples are warped products, for which $\dJ$ can be computed by Lemma~\ref{lem:warped}.

\begin{proposition}\label{prop:warped}
  In the setting of Lemma~\ref{lem:warped}, suppose that $a_- \leq a(t) \leq a_+$ for all $t \in \R$, where $0 < a_- \leq a_+$.
  Then $(M, g)$ is timelike geodesically complete.
  Moreover, the map $\varphi \colon \CM \to \R$, $\varphi(\{t = f(x)\}) \coloneqq (\max f + \min f)/2$, is a $(1, a_+D)$-quasi-isometry,
  and $(\CM, \dJ)$ is $3a_+D$-hyperbolic, where $D \coloneqq \diam(S, h)$.
\end{proposition}

\begin{proof}
  We show timelike geodesic completeness.
  The vertical timelike geodesics $t \mapsto (t, x)$ are complete.
  Since $(S, h)$ is complete and $a \geq a_-$, we have
  $\int_0^\infty (a^2 / (1 + a^2))^{1/2} \, dt = \int_{-\infty}^0 (a^2 / (1 + a^2))^{1/2} \, dt = \infty$,
  so that every other timelike geodesic is complete by a result of Powell~\cite[p.~112]{BEE},
  applied to the warping function $a^2$.

  The proofs of Proposition~\ref{prop:product} (iii) and (iv) use only the following properties:
  every $T \in \CM$ is the graph $\{t = f(x)\}$ of a function $f$ with $\max f - \min f \leq D$,
  and $\dJ(T, T_c) = \max_S \abs{f - c}$ for all $c \in \R$, so that $\dJ(T_c, T_{c'}) = \abs{c - c'}$.
  Hence it suffices to show these properties with $D$ replaced by $a_+D$.
  The formula for $\dJ(T, T_c)$ is Lemma~\ref{lem:warped} (ii).
  For the bound on $\max f - \min f$, let $\eta(t) \coloneqq \int_0^t a(u)^{-1} \, du$.
  Since $1/a_+ \leq \eta' \leq 1/a_-$, the map $\eta$ is a diffeomorphism of $\R$, and $g = a^2(-d\eta^2 + h)$.
  Thus $(M, g)$ is conformal to the product spacetime $(\R \times S, -d\eta^2 + h)$,
  and has the same Cauchy hypersurfaces.
  By Proposition~\ref{prop:product} (i), applied to this product spacetime,
  we have $\abs{\eta(f(x)) - \eta(f(y))} \leq \operatorname{dist}_h(x, y)$ for $T = \{t = f(x)\} \in \CM$.
  Since $\abs{t_1 - t_2} \leq a_+ \abs{\eta(t_1) - \eta(t_2)}$, we obtain $\abs{f(x) - f(y)} \leq a_+ \operatorname{dist}_h(x, y)$,
  and in particular $\max f - \min f \leq a_+D$.
\end{proof}

\begin{example}\label{ex:warped}
  Let $S = \R^n / \Z^n$ be the flat torus and $a(t) \coloneqq 2 + \sin t$.
  By Lemma~\ref{lem:warped} and Proposition~\ref{prop:warped}, $(M, g)$ is a globally hyperbolic, spatially compact, timelike geodesically complete spacetime,
  and $(\CM, \dJ)$ is Gromov hyperbolic and $(1, \varepsilon)$-quasi-isometric to $\R$ for some $\varepsilon \geq 0$.
  The timelike convergence condition fails, since $\Ric(\partial_t, \partial_t) = -n a''/a = n \sin t / (2 + \sin t)$
  is negative for $t \in (\pi, 2\pi)$.

  We show that $(M, g)$ does not split.
  The scalar curvature of $g$ is the function
  \begin{equation}
    \mathrm{scal}_g = 2n \frac{a''}{a} + n(n - 1) \Bigl( \frac{a'}{a} \Bigr)^2
  \end{equation}
  of $t$ alone, which equals $-2n/3$ at $t = \pi/2$ and $2n$ at $t = 3\pi/2$.
  Along a timelike curve parametrized by proper time, $\dot{t}^2 = 1 + a(t)^2 \abs{\dot{x}}_h^2 \geq 1$.
  Hence $t$ takes all real values along every complete timelike geodesic, and $\mathrm{scal}_g$ is not constant along it.
  On the other hand, in a product spacetime $(\R \times \tilde S, -d\tau^2 + \tilde h)$,
  every point lies on a vertical complete timelike geodesic along which the scalar curvature is constant.
  Since isometries preserve the scalar curvature and complete timelike geodesics,
  $(M, g)$ is not isometric to a product spacetime.
  Thus the implications (iii) $\Rightarrow$ (i) and (iv) $\Rightarrow$ (i) of Theorem~\ref{thm:main}
  fail without the timelike convergence condition.
  We note that the vertical lines $t \mapsto (t, x)$ are complete timelike lines by the proof of Lemma~\ref{lem:warped} (ii),
  in accordance with Theorem~\ref{thm:hyperbolic};
  thus the implication (iv) $\Rightarrow$ (i) fails only at the application of the splitting theorem.
\end{example}

\begin{example}\label{ex:desitter}
  Let $n \geq 1$, let $\langle \cdot, \cdot \rangle$ be the Lorentzian scalar product of signature $(-, +, \dots, +)$ on $\R^{n+2}_1$,
  let $\mathbf{S}^{n+1}_1 \coloneqq \{p ; \langle p, p \rangle = 1\}$ be the de Sitter spacetime,
  and let $\mathbf{H}^{n+1} \coloneqq \{x ; \langle x, x \rangle = -1, \ x \text{ future-directed}\}$ be the hyperbolic space,
  whose distance is given by $\cosh d_{\mathbf{H}}(x, y) = -\langle x, y \rangle$.
  The map $\Phi(\tau, u) \coloneqq (\sinh \tau, (\cosh \tau) u)$ is an isometry
  from the warped product $(\R \times \mathbf{S}^n, -d\tau^2 + \cosh^2 \tau \, h_{\mathbf{S}^n})$ onto $\mathbf{S}^{n+1}_1$
  (cf.~\cite[Section 5.2]{HE}), where $h_{\mathbf{S}^n}$ is the round metric.
  Moreover, $\mathbf{S}^{n+1}_1$ is geodesically complete and has constant curvature $1$~\cite[Section 5.3 and Remark 5.2 (2)]{BEE},
  so that $\Ric = n g$ and the timelike convergence condition fails.
  By Lemma~\ref{lem:warped}, the slices $\Phi(\{s\} \times \mathbf{S}^n)$ form a geodesic line in $(\CM, \dJ)$,
  so that (ii) of Theorem~\ref{thm:main} holds.

  To see that (iii) and (iv) fail, let $e_0 \coloneqq (1, 0, \dots, 0) \in \mathbf{H}^{n+1}$.
  For $x \in \mathbf{H}^{n+1}$, choose a linear isometry $\Psi$ of $\R^{n+2}_1$ with $\Psi$$ e_0 = x$.
  Then $\Phi_x \coloneqq \Psi \circ \Phi$ is also an isometry from the above warped product onto $\mathbf{S}^{n+1}_1$,
  and its slices are $\Sigma_{x, s} \coloneqq \Phi_x(\{s\} \times \mathbf{S}^n) = \{p ; \langle p, x \rangle = -\sinh s\}$,
  since $\langle \Phi_x(\tau, u), x \rangle = \langle \Phi(\tau, u), e_0 \rangle = -\sinh \tau$.
  The slices $\Sigma_{x, s}$ form only a special family of Cauchy hypersurfaces of $\mathbf{S}^{n+1}_1$.
  For these slices, we claim that
  \begin{equation}\label{eq:desitter}
    \dJ(\Sigma_{x, s}, \Sigma_{y, t}) = \abs{t - s} + d_{\mathbf{H}}(x, y)
    \qquad \text{for all } x, y \in \mathbf{H}^{n+1} \text{ and } s, t \in \R.
  \end{equation}
  Let $r \coloneqq d_{\mathbf{H}}(x, y)$.
  We identify $\mathbf{S}^{n+1}_1$ with the warped product via $\Phi_x$, under which $\Sigma_{x, s}$ is the slice $\{s\} \times \mathbf{S}^n$.
  Since $\Sigma_{y, t}$ is a Cauchy hypersurface, Lemma~\ref{lem:warped} (i) shows that
  $\Sigma_{y, t} = \{\Phi_x(f(u), u) ; u \in \mathbf{S}^n\}$ for some function $f \colon \mathbf{S}^n \to \R$,
  and Lemma~\ref{lem:warped} (ii) gives $\dJ(\Sigma_{x, s}, \Sigma_{y, t}) = \max_{\mathbf{S}^n} \abs{f - s}$.
  Hence it suffices to show that $\min f = t - r$ and $\max f = t + r$,
  because then $\max_{\mathbf{S}^n} \abs{f - s} = \max\{\abs{t - r - s}, \abs{t + r - s}\} = \abs{t - s} + r$.
  Since $\langle \Phi_x(\tau, u), x \rangle = -\sinh \tau$,
  we have $f(u) = \operatorname{arcsinh}(-\langle \Phi_x(f(u), u), x \rangle)$ for $u \in \mathbf{S}^n$.
  Since $\operatorname{arcsinh}$ is increasing and $\Phi_x(f(u), u)$ ranges over $\Sigma_{y, t}$,
  it is enough to show that $-\langle p, x \rangle$ takes the minimum $\sinh(t - r)$ and the maximum $\sinh(t + r)$ on $\Sigma_{y, t}$.

  Since $\langle x, y \rangle = -\cosh r$ and $\langle x, x \rangle = -1$,
  we can write $x = (\cosh r) y + x'$ with $x' \in y^\perp$ and $\langle x', x' \rangle = \sinh^2 r$.
  Similarly, since $\langle p, y \rangle = -\sinh t$ and $\langle p, p \rangle = 1$ for $p \in \Sigma_{y, t}$, we have
  \begin{equation}
    \Sigma_{y, t} = \{(\sinh t) y + w ; w \in y^\perp, \ \langle w, w \rangle = \cosh^2 t\}.
  \end{equation}
  For $p = (\sinh t) y + w \in \Sigma_{y, t}$, we obtain $-\langle p, x \rangle = \sinh t \cosh r - \langle w, x' \rangle$.
  Since $\langle \cdot, \cdot \rangle$ is positive definite on $y^\perp$,
  the Cauchy--Schwarz inequality shows that $\langle w, x' \rangle$ takes the minimum $-\cosh t \sinh r$
  and the maximum $\cosh t \sinh r$ on $\Sigma_{y, t}$.
  Hence $-\langle p, x \rangle$ takes the minimum $\sinh(t - r)$ and the maximum $\sinh(t + r)$ on $\Sigma_{y, t}$,
  which proves~\eqref{eq:desitter}.

  By~\eqref{eq:desitter}, $(x, s) \mapsto \Sigma_{x, s}$ is an isometric embedding of
  $\mathbf{H}^{n+1} \times \R$ with the $\ell^1$-product distance into $(\CM, \dJ)$.
  Restricting to a geodesic line of $\mathbf{H}^{n+1}$,
  we obtain an isometric copy of $(\R^2, \norm{\cdot}_1)$ in $(\CM, \dJ)$.
  The identity map of $\R^2$ is bi-Lipschitz from $\norm{\cdot}_1$ to the Euclidean norm,
  since $\norm{v}_2 \leq \norm{v}_1 \leq \sqrt{2} \norm{v}_2$.
  Since the Euclidean plane is not Gromov hyperbolic,
  neither is $(\R^2, \norm{\cdot}_1)$ by~\cite[Chapter III.H, Theorem 1.9 and Proposition 1.22]{BH},
  and hence neither is $(\CM, \dJ)$, because the four-point condition passes to subsets.
  Since $(\CM, \dJ)$ is geodesic (Theorem~\ref{thm:injective}), it is not quasi-isometric to $\R$ either.
  Thus the implications (ii) $\Rightarrow$ (iii) and (ii) $\Rightarrow$ (iv) of Theorem~\ref{thm:main}
  fail without the timelike convergence condition.
\end{example}

\begin{example}\label{ex:cylinder}
  Let $M = \R \times \mathbf{S}^1$, where $\mathbf{S}^1 = \R / 2\pi\Z$, and $g = -dt^2 + a(t)^2 d\theta^2$ with $a(t) \coloneqq \tfrac{2}{3} \cosh t$.
  By Lemma~\ref{lem:warped}, $(M, g)$ is globally hyperbolic and spatially compact,
  and it is timelike geodesically complete by the argument of Proposition~\ref{prop:warped}, since $a \geq 2/3$.
  The timelike convergence condition fails, since $\Ric(\partial_t, \partial_t) = -a''/a = -1$.
  With $\eta(t) \coloneqq \int_0^t a(s)^{-1} \, ds$, the metric is conformal to $-d\eta^2 + d\theta^2$ on $(-L/2, L/2) \times \mathbf{S}^1$,
  where $L = \tfrac{3}{2} \int_{\R} \operatorname{sech} s \, ds = \tfrac{3}{2}\pi$,
  and $(\eta, \theta) \ll (\eta', \theta')$ if and only if $\eta' - \eta > \operatorname{dist}(\theta, \theta')$.
  Note that $\pi = \diam \mathbf{S}^1 < L < 2\pi$.

  First, $\CM$ has two ends.
  Since every $T = \{t = f(\theta)\} \in \CM$ is achronal, the above description of $\ll$ shows that $T$ satisfies $\abs{\eta(f(\theta)) - \eta(f(\theta'))} \leq \operatorname{dist}(\theta, \theta') \leq \pi$.
  Choose $t_1 < t_2$ with $\eta(t_2) - \eta(t_1) > \pi$.
  Then no Cauchy hypersurface meets both $\{t < t_1\}$ and $\{t > t_2\}$,
  and $\CM$ has two ends by Proposition~\ref{prop:two-ends}, applied with $\mathcal{T} = t$.

  Second, $N(p)$ is non-compact for every $p = (t_0, \theta_0) \in M$.
  We may assume $\eta(t_0) \geq 0$; the other case is time dual.
  For $t_k \to \infty$, the points $q_k \coloneqq (t_k, \theta_0 + \pi)$ satisfy $0 < \eta(t_k) - \eta(t_0) < L/2 < \pi = \operatorname{dist}(\theta_0, \theta_0 + \pi)$,
  so that $q_k \in N(p)$, and $(q_k)_k$ leaves every compact set, since $t(q_k) = t_k \to \infty$.
  Thus the point condition of~\cite{Ba} fails, although $\CM$ has two ends.
\end{example}

\section{Local geodesics and the failure of convexity}\label{sec:convexity}

In this section, we show that distance curves of smooth spacelike Cauchy hypersurfaces are local geodesics
in $(\CM, \dJ)$, and deduce that $(\CM, \dJ)$ is neither CAT(0) nor Busemann convex in general,
and never under the timelike convergence condition,
although it is injective by Theorem~\ref{thm:injective}.
We then compare $\dJ$ with the $L^2$-metric of~\cite{Mon}.

\begin{lemma}\label{lem:local-geodesic}
  Let $\Sigma$ be a smooth, spacelike, acausal compact Cauchy hypersurface.
  Then there exists $\varepsilon > 0$ such that $\dJ(\Sigma_s, \Sigma_t) = \abs{t - s}$ for all $s, t \in [-\varepsilon, \varepsilon]$.
  Consequently, the distance curve $c_\Sigma$ is a local geodesic in $(\CM, \dJ)$.
\end{lemma}

\begin{proof}
  For $x \in \Sigma$, let $\gamma_x$ be the timelike geodesic with $\gamma_x(0) = x$ whose initial velocity
  is the future-directed unit normal of $\Sigma$ at $x$, and set $E(x, t) \coloneqq \gamma_x(t)$.
  By~\cite[Proposition 7.26]{On}, $E$ is injective on a neighborhood of $\Sigma \times \{0\}$,
  which contains $\Sigma \times (-\varepsilon_0, \varepsilon_0)$ for some $\varepsilon_0 > 0$ since $\Sigma$ is compact.
  The function $(x, t) \mapsto d(\Sigma, \gamma_x(t))$ is continuous by Fact~\ref{fact:cont},
  and it vanishes for $t = 0$ since $\Sigma$ is achronal; the same holds for $(x, t) \mapsto d(\gamma_x(-t), \Sigma)$.
  By the compactness of $\Sigma$, there exists $\varepsilon \in (0, \varepsilon_0)$ such that
  $d(\Sigma, \gamma_x(t)) < \varepsilon_0$ and $d(\gamma_x(-t), \Sigma) < \varepsilon_0$ for all $x \in \Sigma$ and $t \in [0, \varepsilon]$.

  Let $x \in \Sigma$ and $0 < t \leq \varepsilon$, and set $p \coloneqq \gamma_x(t)$.
  Then $\ell \coloneqq d(\Sigma, p)$ satisfies $t \leq \ell < \varepsilon_0$, since $\gamma_x|_{[0, t]}$ is a causal curve
  of length $t$ from $\Sigma$ to $p$.
  By Lemma~\ref{lem:cauchy} (iii) and Fact~\ref{fact:gh}, there exist $y \in \Sigma$ with $d(y, p) = \ell$
  and a maximizing geodesic from $y$ to $p$ of length $\ell$.
  Since this geodesic maximizes the length among causal curves from $\Sigma$ to $p$,
  it is a timelike geodesic orthogonal to $\Sigma$ by~\cite[Proposition 2.21]{MiLR}.
  Hence it coincides with $\gamma_y|_{[0, \ell]}$, and $E(y, \ell) = p = E(x, t)$.
  The injectivity of $E$ gives $y = x$ and $\ell = t$.
  Thus $\gamma_x(t) \in \Sigma_t$, and similarly $\gamma_x(-t) \in \Sigma_{-t}$, for $0 \leq t \leq \varepsilon$.

  Now let $-\varepsilon \leq s \leq 0 \leq t \leq \varepsilon$. Then
  \begin{equation}
    d(\Sigma_s, \Sigma_t) \geq d(\gamma_x(s), \gamma_x(t)) \geq t - s,
  \end{equation}
  while $d(\Sigma_s, \Sigma_t) \leq t - s$ by Fact~\ref{fact:lengths}.
  By Lemma~\ref{lem:monotone-family}, we obtain $\dJ(\Sigma_s, \Sigma_t) = t - s$.
  Together with Proposition~\ref{prop:parallel} (i), this proves the first claim.
  For $t_0 \neq 0$, the distance curve is isometric on a neighborhood of $t_0$ by Proposition~\ref{prop:parallel} (i),
  so that $c_\Sigma$ is a local geodesic.
\end{proof}

Thus the distance curve of a smooth spacelike Cauchy hypersurface is always a local geodesic,
and the defect $D_\infty(\Sigma)$ arises only globally.
Since local geodesics are geodesics in CAT(0) and Busemann convex spaces (Fact~\ref{fact:cat0}), this leads to the following.

\begin{theorem}\label{thm:convexity}
  Let $(M, g)$ be a timelike geodesically complete, spatially compact, globally hyperbolic spacetime.
  If some smooth, spacelike, acausal compact Cauchy hypersurface $\Sigma$ satisfies $D_\infty(\Sigma) > 0$,
  then $(\CM, \dJ)$ is neither a CAT(0) space nor Busemann convex.
  In particular, for every product spacetime $(\R \times S, -dt^2 + h)$ with $(S, h)$ compact,
  $(\CM, \dJ)$ is a proper geodesic Gromov hyperbolic space which is neither CAT(0) nor Busemann convex.
\end{theorem}

\begin{proof}
  By Proposition~\ref{prop:parallel} (ii), we have $D_\infty(\Sigma) = \lim_{k \to \infty} D(k, k) > 0$,
  so that there exists $L > 0$ with $D(L, L) > 0$.
  By Lemma~\ref{lem:local-geodesic}, the curve $c_\Sigma|_{[-L, L]}$ is a local geodesic.
  However, $\dJ(\Sigma_{-L}, \Sigma_L) = 2L - D(L, L) < 2L$ by Proposition~\ref{prop:parallel} (ii),
  so that $c_\Sigma|_{[-L, L]}$ is not a geodesic.
  If $(\CM, \dJ)$ were a CAT(0) space or Busemann convex, then every local geodesic would be a geodesic by Fact~\ref{fact:cat0},
  which is a contradiction.

  Finally, let $(M, g) = (\R \times S, -dt^2 + h)$ with $(S, h)$ compact; $S$ is connected since $M$ is connected.
  This spacetime satisfies our standing assumptions by the paragraph before Proposition~\ref{prop:product}.
  By Proposition~\ref{prop:product} (iv) and (v), $(\CM, \dJ)$ is a proper geodesic Gromov hyperbolic space.
  Let $f \colon S \to \R$ be a smooth non-constant function with $\abs{df}_h < 1$,
  and let $\Sigma \coloneqq \{t = f(x)\}$, which is a smooth, spacelike, compact hypersurface.
  Since $f$ is $1$-Lipschitz, $\Sigma \in \CM$ by Proposition~\ref{prop:product} (i).
  Along a future-directed causal curve $s \mapsto (t(s), x(s))$, we have $\dot{t} \geq \abs{\dot{x}}_h$ and $\dot{t} > 0$,
  so that $\frac{d}{ds}(t - f(x)) \geq \dot{t} - \abs{df}_h \abs{\dot{x}}_h > 0$, and hence $\Sigma$ is acausal.
  Thus $\Sigma$ is an acausal Cauchy hypersurface,
  and $D_\infty(\Sigma) = \max f - \min f > 0$ by Fact~\ref{fact:product} and Proposition~\ref{prop:parallel} (ii),
  so that the first statement applies.
\end{proof}

\begin{corollary}\label{cor:convexity}
  In the setting of Theorem~\ref{thm:convexity}, suppose in addition that
  $(M, g)$ satisfies the timelike convergence condition.
  Then $(\CM, \dJ)$ is neither a CAT(0) space nor Busemann convex.
\end{corollary}

\begin{proof}
  If $(M, g)$ splits, then it is isometric to a product spacetime $(\R \times S, -dt^2 + h)$ with $(S, h)$ compact.
  An isometry between spacetimes preserves the Lorentzian distance and maps Cauchy hypersurfaces to Cauchy hypersurfaces,
  and hence induces an isometry between the spaces of Cauchy hypersurfaces.
  Thus the claim follows from the last statement of Theorem~\ref{thm:convexity}.
  Otherwise, $\lim_{k \to \infty} \Delta_k = \infty$ by Fact~\ref{fact:defect},
  so that $D_\infty(\Sigma) = \infty$ for every acausal compact Cauchy hypersurface $\Sigma$
  by Proposition~\ref{prop:parallel} (ii).
  Since a smooth, spacelike, acausal compact Cauchy hypersurface exists~\cite{BS}, \cite[Theorem 4.120]{MiLR},
  the first statement of Theorem~\ref{thm:convexity} applies.
\end{proof}

\begin{remark}\label{rem:monclair}
  Monclair~\cite{Mon} showed that the $L^2$-metric on the space of smooth spacelike Cauchy hypersurfaces
  has non-positive sectional curvature, and asked whether this space, or its metric completion, is a CAT(0) space~\cite[Question 1.4]{Mon}.
  Theorem~\ref{thm:convexity} and Corollary~\ref{cor:convexity} show that $\dJ$ behaves differently.
  Since Lange and Peteranderl describe $\dJ$ as an $L^\infty$-Finsler analogue of Monclair's metric~\cite[Section 1]{LP},
  the failure of CAT(0) and Busemann convexity itself may not be surprising.
  The point of Theorem~\ref{thm:convexity} is that this failure is detected by the distance curves, which are local geodesics (Lemma~\ref{lem:local-geodesic}),
  that it is governed by the defect $D_\infty(\Sigma)$,
  and that, by Corollary~\ref{cor:convexity}, it occurs for every spacetime satisfying the timelike convergence condition.
\end{remark}

\section{Concluding remarks}\label{sec:remarks}

By Corollary~\ref{cor:coarse}, Bartnik's conjecture is equivalent
to a statement on the coarse geometry of $(\CM, \dJ)$
for all spacetimes under consideration,
and by Corollary~\ref{cor:type}, even to a statement on the topology of $\CM$.
We do not know whether the timelike convergence condition
implies a coarse property of $(\CM, \dJ)$ of this kind.
Such a property cannot come from CAT(0) or Busemann convexity,
since under the timelike convergence condition $(\CM, \dJ)$
is neither a CAT(0) space nor Busemann convex
(Corollary~\ref{cor:convexity}).
On the other hand, $(\CM, \dJ)$ is injective (Theorem~\ref{thm:injective}),
a property shared with non-positively curved spaces, but this holds without any curvature condition.
Thus the timelike convergence condition would have to be reflected
in $(\CM, \dJ)$ through a property of a different nature.
We do not know the homeomorphism type of $\CM$ in a hypothetical counterexample to Bartnik's conjecture.
It has exactly one end by Corollary~\ref{cor:ends}, and a natural question is whether it is homeomorphic to $[0, \infty) \times Q$.

The metric $\dJ$ does not rely on smoothness.
In~\cite[Section 5]{LP}, it is defined on the space of
Cauchy sets of a Lorentzian pre-length space,
and variants of the results of~\cite{LP} are proved
in the synthetic settings of
Kunzinger--S\"amann~\cite{KS18}, Braun--McCann~\cite{BM26} and Bykov--Minguzzi--Suhr~\cite{BMS25}. 
In the present article, smoothness is used in the following places:
the graph representation in the proof of Proposition~\ref{prop:proper},
as in~\cite[Proposition 5.23]{LP};
the results of~\cite{Mi} on parallel sets, which are formulated for smooth spacetimes;
the extension of maximizing segments to complete timelike geodesics
in the proof of Theorem~\ref{thm:quasi};
and the normal exponential map in Lemma~\ref{lem:local-geodesic}.
On the other hand, the lattice structure of $\CM$ (Lemma~\ref{lem:lattice}) is obtained by causal arguments,
and the fact that parallel sets are Cauchy sets has a synthetic counterpart~\cite[Lemma 5.16]{LP}.
Moreover, the only curvature input in Theorem~\ref{thm:main}
is the Lorentzian splitting theorem,
which has been extended to globally hyperbolic Lorentzian length spaces
with non-negative timelike curvature~\cite{BORS}
and to continuously differentiable metrics~\cite{BGMOS}.
The implication (iii) $\Rightarrow$ (i) of Corollary~\ref{cor:type} consists of two steps.
The first step, from the compactness of $\CM(p)$ to two ends, is a topological argument based on Lemma~\ref{lem:topology}.
The second step, from two ends to the splitting, goes through Theorem~\ref{thm:main},
and hence uses the smooth structure as listed above, in particular in the proof of Theorem~\ref{thm:quasi}.
It might be worth investigating whether Theorems~\ref{thm:injective}, \ref{thm:quasi} and~\ref{thm:main}
have counterparts in such settings.

\begin{acknowledgements}
The author would like to thank Christian Lange and Jonas W.~Peteranderl
for stimulating conversations on the space of Cauchy hypersurfaces
during the conference ``Curvature, Optimal Transport, and Gravity''
at the University of Vienna in September 2026.
This research was supported in part by the Austrian Science Fund (FWF)
[Grants DOI 10.55776/EFP6, 10.55776/STA32].
\end{acknowledgements}

\begin{aistatement}
Claude and ChatGPT were used interactively in the preparation of this manuscript.
The author independently reviewed and revised all mathematical statements, proofs and text,
and takes full responsibility for the final content.
\end{aistatement}

\bibliographystyle{alpha}
\bibliography{bibliography}

\end{document}